\documentclass[12pt]{article}
\usepackage[left=1in,right=1in,top=1in,bottom=1in]{geometry}
\usepackage[utf8]{inputenc}
\usepackage[T1]{fontenc}
\usepackage[colorinlistoftodos]{todonotes}

\usepackage{amssymb,amsmath,amsthm,comment,enumerate,hyperref,multicol,mathdots,tikz,ytableau,mathtools}
\usepackage{mlmodern}
\usetikzlibrary{cd}

\theoremstyle{definition}
\newtheorem{theorem}{Theorem}[section]
\newtheorem{lemma}[theorem]{Lemma}
\newtheorem{proposition}[theorem]{Proposition}
\newtheorem{corollary}[theorem]{Corollary}
\newtheorem{definition}[theorem]{Definition}

\newtheorem{problem}[theorem]{Problem}

\newcommand{\Z}{\mathbb{Z}}

\newcommand{\BPD}[1]{\text{BPD}(#1)}

\definecolor{UF2}{HTML}{FA4616}
\definecolor{UF}{HTML}{0021A5}
\newcommand{\cross}[1][UF]{\begin{tikzpicture}[scale = .4]
	\draw[lightgray] (0,0) -- (1,0) -- (1,1) -- (0,1) -- (0,0);
	\draw[thick,#1!100!black] (.5,0) -- (.5,1);
	\draw[thick,#1!100!black] (0,.5) -- (1,.5);
\end{tikzpicture}}
\newcommand{\hwire}[1][UF]{\begin{tikzpicture}[scale = .4]
	\draw[lightgray] (0,0) -- (1,0) -- (1,1) -- (0,1) -- (0,0);
	\draw[thick,#1!100!black] (0,.5) -- (1,.5);
\end{tikzpicture}}

\newcommand{\hewire}[1][UF]{\begin{tikzpicture}[scale = .4]
	\draw[lightgray] (0,0) -- (1,0) -- (1,1) -- (0,1) -- (0,0);
	\draw[thick,dotted, #1!100!black] (0,.5) -- (1,.5);
    \draw[thick,#1!100!black] (0,0.5) to[out=0,in=-90] (0.5,1);
\end{tikzpicture}}

\newcommand{\vwire}[1][UF]{\begin{tikzpicture}[scale = .4]
	\draw[lightgray] (0,0) -- (1,0) -- (1,1) -- (0,1) -- (0,0);
	\draw[thick,#1!100!black] (.5,0) -- (.5,1);
\end{tikzpicture}}

\newcommand{\vewire}[1][UF]{\begin{tikzpicture}[scale = .4]
	\draw[lightgray] (0,0) -- (1,0) -- (1,1) -- (0,1) -- (0,0);
	\draw[thick,dotted,#1!100!black] (.5,0) -- (.5,1);
    \draw[thick,#1!100!black] (0,0.5) to[out=0,in=-90] (0.5,1);

\end{tikzpicture}}

\newcommand{\nowire}{\begin{tikzpicture}[scale = .4]
	\draw[lightgray] (0,0) -- (1,0) -- (1,1) -- (0,1) -- (0,0);
\end{tikzpicture}}

\newcommand{\jay}[1][UF]{\begin{tikzpicture}[scale = .4]
	\draw[lightgray] (0,0) -- (1,0) -- (1,1) -- (0,1) -- (0,0);
	\draw[thick,#1!100!black] (0,0.5) to[out=0,in=-90] (0.5,1);
\end{tikzpicture}}

\newcommand{\are}[1][UF]{\begin{tikzpicture}[scale = .4]
	\draw[lightgray] (0,0) -- (1,0) -- (1,1) -- (0,1) -- (0,0);
	\draw[thick,#1!100!black] (0.5,0) to[out=90,in=180] (1,0.5);
\end{tikzpicture}}
\newcommand{\inlinecross}{\begin{tikzpicture}[scale = .3]
	\draw[lightgray] (0,0) -- (1,0) -- (1,1) -- (0,1) -- (0,0);
	\draw[thick,UF!100!black] (.5,0) -- (.5,1);
	\draw[thick,UF!100!black] (0,.5) -- (1,.5);
\end{tikzpicture}}
\newcommand{\inlinenowire}{\begin{tikzpicture}[scale = .3]
	\draw[lightgray] (0,0) -- (1,0) -- (1,1) -- (0,1) -- (0,0);
\end{tikzpicture}}
\newcommand{\inlinebump}{\begin{tikzpicture}[scale = .3]
	\draw[lightgray] (0,0) -- (1,0) -- (1,1) -- (0,1) -- (0,0);
	\draw[thick,UF!100!black] (0.5,0) to[out=90,in=180] (1,0.5);
	\draw[thick,UF!100!black] (0,0.5) to[out=0,in=-90] (0.5,1);
\end{tikzpicture}}
\newcommand{\inlineare}{\begin{tikzpicture}[scale = .3]
	\draw[lightgray] (0,0) -- (1,0) -- (1,1) -- (0,1) -- (0,0);
	\draw[thick,UF!100!black] (0.5,0) to[out=90,in=180] (1,0.5);
\end{tikzpicture}}
\newcommand{\inlinejay}{\begin{tikzpicture}[scale = .3]
	\draw[lightgray] (0,0) -- (1,0) -- (1,1) -- (0,1) -- (0,0);
	\draw[thick,UF!100!black] (0,0.5) to[out=0,in=-90] (0.5,1);
\end{tikzpicture}}
\newcommand{\inlinehwire}{\begin{tikzpicture}[scale = .3]
	\draw[lightgray] (0,0) -- (1,0) -- (1,1) -- (0,1) -- (0,0);
	\draw[thick,UF!100!black] (0,.5) -- (1,.5);
\end{tikzpicture}}
\newcommand{\inlinevwire}{\begin{tikzpicture}[scale = .3]
	\draw[lightgray] (0,0) -- (1,0) -- (1,1) -- (0,1) -- (0,0);
	\draw[thick,UF!100!black] (.5,0) -- (.5,1);
\end{tikzpicture}}

\begin{document}
\title{A Graph on Reduced Bumpless Pipe Dreams}
\author{
Joshua Arroyo
\footnote{Georgia State University; Atlanta, GA} 
\and
Adam Gregory
\footnote{Western Carolina University; Cullowhee, NC 
} 
\and 
Prashreet Poudel
\footnotemark[2]
}
\date{\today}
\maketitle

\begin{abstract}
Reduced bumpless pipe dreams are combinatorial objects used to study Schubert polynomials. 
They resemble tilings of a square grid using one of six pieces subject to a few constraints, and they can be partitioned into a collection which is indexed by permutations~$w$.
The set of reduced bumpless pipe dreams for $w$ can be generated via certain local moves that naturally give rise to a graph structure we associate with $w$.
In this note, we provide a necessary pattern-avoidance condition on the permutations whose corresponding reduced bumpless pipe dream graph is a tree.
\end{abstract}


\section{Introduction}
\label{s:intro}
A {\em reduced bumpless pipe dream of size} $n$ is a tiling of an $n \times n$ grid using tiles from
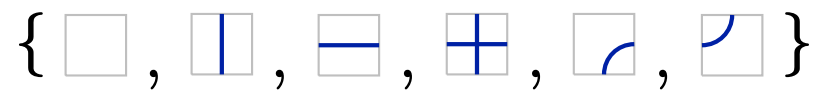
\[
\raisebox{3pt}{\{} \; \nowire \; , \; \vwire \; , \; \hwire \; , \; \cross \; , \; \are \; , \; \jay \; \raisebox{3pt}{\}}
\]
such that $n$ pipes enter from the bottom edge and exit to the right edge, with no two pipes crossing each other more than once. Satisfying this last condition is what the word ``reduced'' is implying. Lam, Lee, and Shimozono introduced reduced bumpless pipe dreams as a combinatorial model for studying Schubert polynomials~\cite{lam2021back}. 

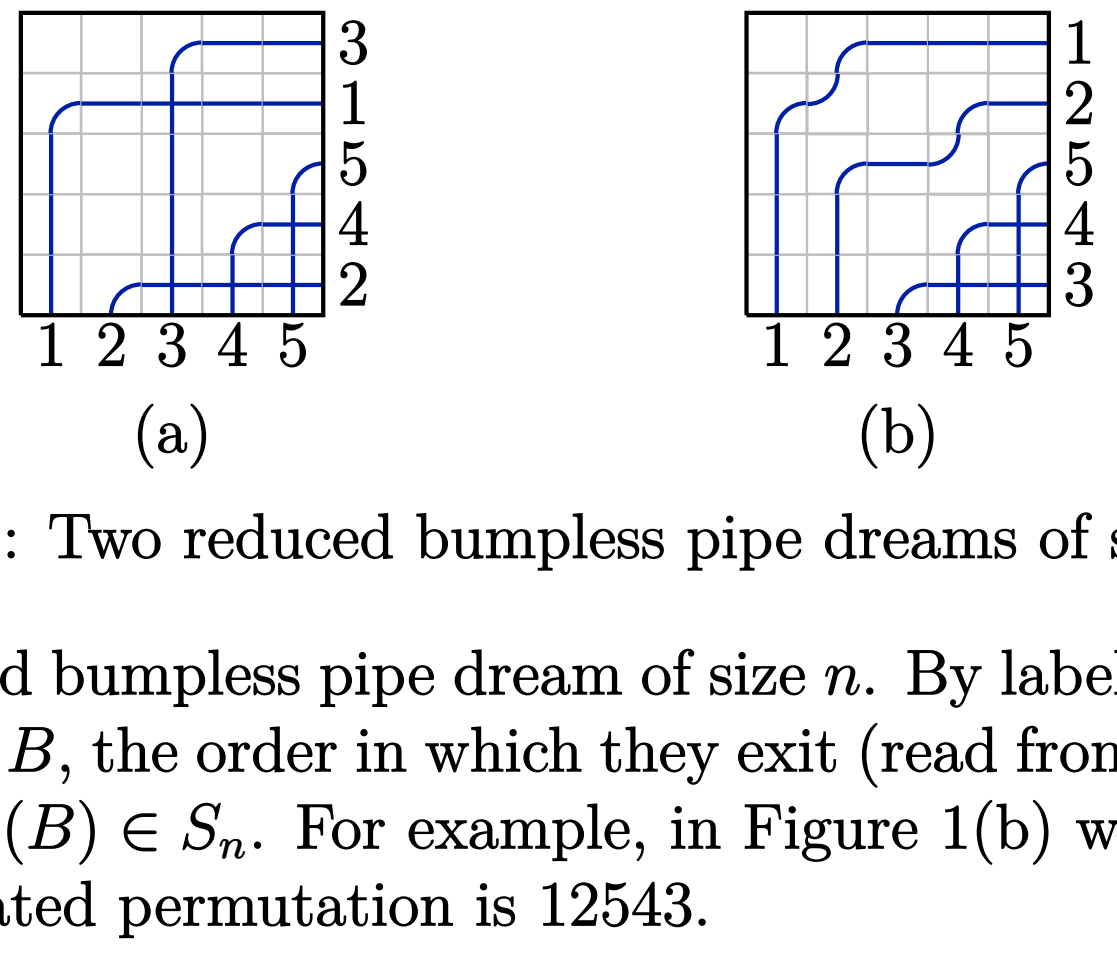
\begin{figure}[h!]
\[
\begin{tikzpicture}[scale=0.4]
\node at (0,0) {\vwire};
\node at (1,0) {\are};
\node at (2,0) {\cross};
\node at (3,0) {\cross};
\node at (4,0) {\cross};

\node at (0,1) {\vwire};
\node at (1,1) {\nowire};
\node at (2,1) {\vwire};
\node at (3,1) {\are};
\node at (4,1) {\cross};

\node at (0,2) {\vwire};
\node at (1,2) {\nowire};
\node at (2,2) {\vwire};
\node at (3,2) {\nowire};
\node at (4,2) {\are};

\node at (0,3) {\are};
\node at (1,3) {\hwire};
\node at (2,3) {\cross};
\node at (3,3) {\hwire};
\node at (4,3) {\hwire};

\node at (0,4) {\nowire};
\node at (1,4) {\nowire};
\node at (2,4) {\are};
\node at (3,4) {\hwire};
\node at (4,4) {\hwire};
\draw[thick] (-0.5,-0.5) -- (-0.5,4.5) -- (4.5,4.5) -- (4.5,-0.5) -- (-0.5,-0.5);

\node at (0,-1) {1};
\node at (1,-1) {2};
\node at (2,-1) {3};
\node at (3,-1) {4};
\node at (4,-1) {5};

\node at (5,4) {3};
\node at (5,3) {1};
\node at (5,2) {5};
\node at (5,1) {4};
\node at (5,0) {2};

\node at (2,-2.5) {(a)};
\node at (14,-2.5) {(b)};


\node at (12,0) {\vwire};
\node at (13,0) {\vwire};
\node at (14,0) {\are};
\node at (15,0) {\cross};
\node at (16,0) {\cross};

\node at (12,1) {\vwire};
\node at (13,1) {\vwire};
\node at (14,1) {\nowire};
\node at (15,1) {\are};
\node at (16,1) {\cross};

\node at (12,2) {\vwire};
\node at (13,2) {\are};
\node at (14,2) {\hwire};
\node at (15,2) {\jay};
\node at (16,2) {\are};

\node at (12,3) {\are};
\node at (13,3) {\jay};
\node at (14,3) {\nowire};
\node at (15,3) {\are};
\node at (16,3) {\hwire};

\node at (12,4) {\nowire};
\node at (13,4) {\are};
\node at (14,4) {\hwire};
\node at (15,4) {\hwire};
\node at (16,4) {\hwire};
\draw[thick] (11.5,-0.5) -- (11.5,4.5) -- (16.5,4.5) -- (16.5,-0.5) -- (11.5,-0.5);

\node at (12,-1) {1};
\node at (13,-1) {2};
\node at (14,-1) {3};
\node at (15,-1) {4};
\node at (16,-1) {5};

\node at (17,4) {1};
\node at (17,3) {2};
\node at (17,2) {5};
\node at (17,1) {4};
\node at (17,0) {3};

\end{tikzpicture} 
\]
\vspace{-0.45in}
\caption{Two reduced bumpless pipe dreams of size five.}
\label{fig:two-BPDs}
\end{figure}

Suppose $B$ is a reduced bumpless pipe dream of size $n$. By labeling the pipes $1$-to-$n$ from left-to-right as they enter $B$, the order in which they exit (read from top-to-bottom) gives its associated permutation $\pi(B) \in S_n$. For example, in Figure~\ref{fig:two-BPDs}(b) we see a reduced bumpless pipe dream whose associated permutation is 12543. 

\newpage

For a fixed $n \geq 1$, let $\text{BPD}(n)$ denote the set of reduced bumpless pipe dreams of size $n$. Now, for a fixed $w \in S_n$, let 
\[
\text{BPD}(w) = \{B \in \text{BPD}(n) : \pi(B) = w\}.
\]
Note that $|\text{BPD}(w)| \geq 1$ for any permutation $w$. Indeed, with $w$ fixed, pipe $i$ can enter from the bottom edge and travel vertically using $\inlinevwire$ tiles until reaching the row in which $w$ has prescribed that $i$ must exit. Then, pipe $i$ can turn right using a $\inlineare$ tile and continue traveling right using $\inlinehwire$ tiles until exiting. Doing this for each $i \in \{1,\ldots,n\}$ is guaranteed to build a reduced bumpless pipe dream $\overline{B}_w$ whose associated permutation is $w$. For example, in Figure~\ref{fig:two-BPDs}(a) we see $\overline{B}_{31542}$. Note that $\overline{B}_w$ is usually called the {\em Rothe bumpless pipe dream} for $w$ because of the connection between the locations of its $\inlinenowire$ tiles and what is known as the ``Rothe diagram'' of $w$. No $\inlinejay$ tiles are used in $\overline{B}_w$ for any $w$.  \\

In~\cite[Proposition 5.3]{lam2021back} it was shown that each $B \in \text{BPD}(w)$ can be obtained from $\overline{B}_w$ by a sequence of what are called ``droop moves.'' Informally, a {\em droop move} on a reduced bumpless pipe dream is a reassignment of the tiles matching one of the following forms:
    
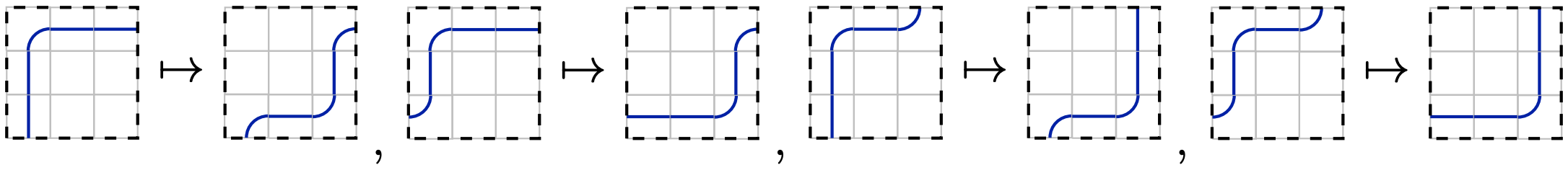
\[
\begin{tikzpicture}[scale= .4]
\node at (0,0) {\vwire};
\node at (1,0) {\nowire};
\node at (2,0) {\nowire};

\node at (0,1) {\vwire};
\node at (1,1) {\nowire};
\node at (2,1) {\nowire};

\node at (0,2) {\are};
\node at (1,2) {\hwire};
\node at (2,2) {\hwire};

\draw[thick, dashed] (-0.5,-0.5) -- (-0.5,2.5) -- (2.5,2.5) -- (2.5,-0.5) -- (-0.5,-0.5);

\node at (5,0) {\are};
\node at (6,0) {\hwire};
\node at (7,0) {\jay};

\node at (5,1) {\nowire};
\node at (6,1) {\nowire};
\node at (7,1) {\vwire};

\node at (5,2) {\nowire};
\node at (6,2) {\nowire};
\node at (7,2) {\are};

\node at (3.5, 1) {$\mapsto$};

\draw[thick, dashed] (4.5,-0.5) -- (4.5,2.5) -- (7.5,2.5) -- (7.5,-0.5) -- (4.5,-0.5);

\end{tikzpicture},
\begin{tikzpicture}[scale= .4]
\node at (0,0) {\jay};
\node at (1,0) {\nowire};
\node at (2,0) {\nowire};

\node at (0,1) {\vwire};
\node at (1,1) {\nowire};
\node at (2,1) {\nowire};

\node at (0,2) {\are};
\node at (1,2) {\hwire};
\node at (2,2) {\hwire};

\draw[thick, dashed] (-0.5,-0.5) -- (-0.5,2.5) -- (2.5,2.5) -- (2.5,-0.5) -- (-0.5,-0.5);

\node at (5,0) {\hwire};
\node at (6,0) {\hwire};
\node at (7,0) {\jay};

\node at (5,1) {\nowire};
\node at (6,1) {\nowire};
\node at (7,1) {\vwire};

\node at (5,2) {\nowire};
\node at (6,2) {\nowire};
\node at (7,2) {\are};

\node at (3.5, 1) {$\mapsto$};

\draw[thick, dashed] (4.5,-0.5) -- (4.5,2.5) -- (7.5,2.5) -- (7.5,-0.5) -- (4.5,-0.5);
\end{tikzpicture},
\begin{tikzpicture}[scale= .4]
\node at (0,0) {\vwire};
\node at (1,0) {\nowire};
\node at (2,0) {\nowire};

\node at (0,1) {\vwire};
\node at (1,1) {\nowire};
\node at (2,1) {\nowire};

\node at (0,2) {\are};
\node at (1,2) {\hwire};
\node at (2,2) {\jay};

\draw[thick, dashed] (-0.5,-0.5) -- (-0.5,2.5) -- (2.5,2.5) -- (2.5,-0.5) -- (-0.5,-0.5);

\node at (5,0) {\are};
\node at (6,0) {\hwire};
\node at (7,0) {\jay};

\node at (5,1) {\nowire};
\node at (6,1) {\nowire};
\node at (7,1) {\vwire};

\node at (5,2) {\nowire};
\node at (6,2) {\nowire};
\node at (7,2) {\vwire};

\node at (3.5, 1) {$\mapsto$};

\draw[thick, dashed] (4.5,-0.5) -- (4.5,2.5) -- (7.5,2.5) -- (7.5,-0.5) -- (4.5,-0.5);

\end{tikzpicture},
\begin{tikzpicture}[scale= .4]
\node at (0,0) {\jay};
\node at (1,0) {\nowire};
\node at (2,0) {\nowire};

\node at (0,1) {\vwire};
\node at (1,1) {\nowire};
\node at (2,1) {\nowire};

\node at (0,2) {\are};
\node at (1,2) {\hwire};
\node at (2,2) {\jay};

\draw[thick, dashed] (-0.5,-0.5) -- (-0.5,2.5) -- (2.5,2.5) -- (2.5,-0.5) -- (-0.5,-0.5);

\node at (5,0) {\hwire};
\node at (6,0) {\hwire};
\node at (7,0) {\jay};

\node at (5,1) {\nowire};
\node at (6,1) {\nowire};
\node at (7,1) {\vwire};

\node at (5,2) {\nowire};
\node at (6,2) {\nowire};
\node at (7,2) {\vwire};

\node at (3.5, 1) {$\mapsto$};

\draw[thick, dashed] (4.5,-0.5) -- (4.5,2.5) -- (7.5,2.5) -- (7.5,-0.5) -- (4.5,-0.5);

\end{tikzpicture}
\]
where unpictured pipes may pass horizontally and vertically through the region, but no $\inlineare$ or $\inlinejay$ tiles appear within the region. Such a move does not affect the order in which the pipes exit nor cause two pipes to cross more than once. Hence, a droop move on $B \in \text{BPD}(w)$ results in a new reduced bumpless pipe dream $B' \in \text{BPD}(w)$. A droop move is called {\em minimal} (or a {\em min-droop} for short) if there are no $\inlinenowire$ tiles in the drooping region other than the bottom right corner. Note that these are slightly different than the ``min-droops'' of~\cite[Definition 3.1]{huang2022bumpless}. \\

We introduce the {\em reduced bumpless pipe dream graph} for $w$, denoted here by $H_w$, to be the graph whose vertices are $\text{BPD}(w)$, with an edge between any $B_1,B_2 \in \text{BPD}(w)$ whenever $B_1$ and $B_2$ are related via a min-droop move. We show (Lemma~\ref{l:droop-imply-min-droop}) that any droop move is a sequence of min-droops, hence these $H_w$ graphs are connected for any $w$. \\

\begin{figure}[h!]
\centering
\[
\begin{tikzpicture}[scale =.4]\node at (0,3) {\are};\node at (1,3) {\hwire};\node at (2,3) {\hwire};\node at (3,3) {\hwire};\node at (0,2) {\vwire};\node at (1,2) {\nowire};\node at (2,2) {\nowire};\node at (3,2) {\are};\node at (0,1) {\vwire};\node at (1,1) {\nowire};\node at (2,1) {\are};\node at (3,1) {\cross};\node at (0,0) {\vwire};\node at (1,0) {\are};\node at (2,0) {\cross};\node at (3,0) {\cross};

\node at (8,3) {\nowire};\node at (9,3) {\are};\node at (10,3) {\hwire};\node at (11,3) {\hwire};\node at (8,2) {\are};\node at (9,2) {\jay};\node at (10,2) {\nowire};\node at (11,2) {\are};\node at (8,1) {\vwire};\node at (9,1) {\nowire};\node at (10,1) {\are};\node at (11,1) {\cross};\node at (8,0) {\vwire};\node at (9,0) {\are};\node at (10,0) {\cross};\node at (11,0) {\cross};

\node at (16,7) {\nowire};\node at (17,7) {\nowire};\node at (18,7) {\are};\node at (19,7) {\hwire};\node at (16,6) {\are};\node at (17,6) {\hwire};\node at (18,6) {\jay};\node at (19,6) {\are};\node at (16,5) {\vwire};\node at (17,5) {\nowire};\node at (18,5) {\are};\node at (19,5) {\cross};\node at (16,4) {\vwire};\node at (17,4) {\are};\node at (18,4) {\cross};\node at (19,4) {\cross};

\node at (16,-1) {\nowire};\node at (17,-1) {\are};\node at (18,-1) {\hwire};\node at (19,-1) {\hwire};\node at (16,-2) {\nowire};\node at (17,-2) {\vwire};\node at (18,-2) {\nowire};\node at (19,-2) {\are};\node at (16,-3) {\are};\node at (17,-3) {\jay};\node at (18,-3) {\are};\node at (19,-3) {\cross};\node at (16,-4) {\vwire};\node at (17,-4) {\are};\node at (18,-4) {\cross};\node at (19,-4) {\cross};

\node at (24,3) {\nowire};\node at (25,3) {\nowire};\node at (26,3) {\are};\node at (27,3) {\hwire};\node at (24,2) {\nowire};\node at (25,2) {\are};\node at (26,2) {\jay};\node at (27,2) {\are};\node at (24,1) {\are};\node at (25,1) {\jay};\node at (26,1) {\are};\node at (27,1) {\cross};\node at (24,0) {\vwire};\node at (25,0) {\are};\node at (26,0) {\cross};\node at (27,0) {\cross};

\draw[thick] (4,1.5) -- (7,1.5);
\draw[thick] (12,2.5) -- (15,5.5);
\draw[thick] (12,0.5) -- (15,-2.5);
\draw[thick] (20, 5.5) -- (23,2.5);
\draw[thick] (20,-2.5) -- (23,0.5);
\end{tikzpicture}
\]
\caption{The reduced bumpless pipe dream graph $H_{1432}$.}
\label{fig:H1432}
\end{figure}
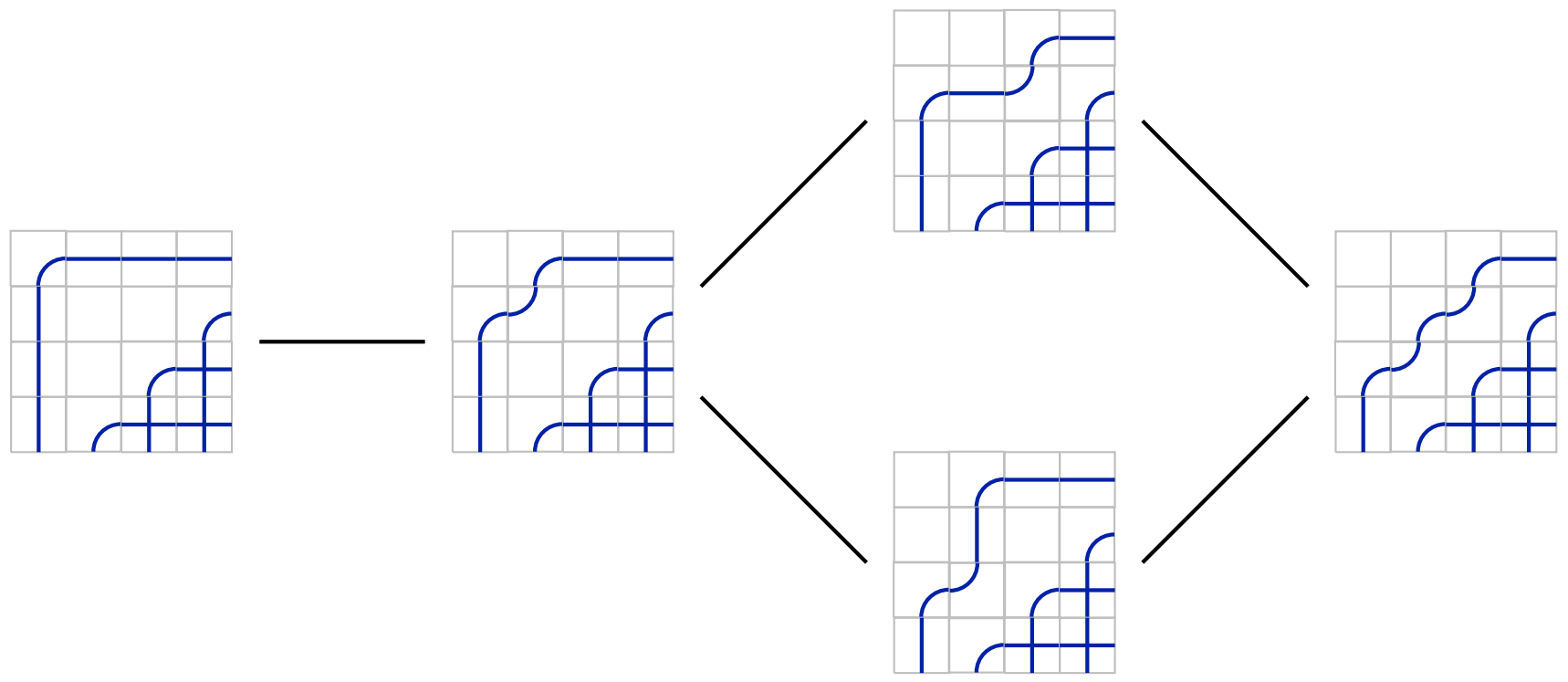


Our main result is the following pattern-avoidance condition on $w$ that is necessary for the reduced bumpless pipe dream graph $H_w$ to be a tree. Recall that $w \in S_n$ {\em avoids} $u \in S_m$ if there are no indices $i_1 < \ldots < i_m$ with the images $w(i_1), \ldots, w(i_m)$ appearing in the same relative order as in $u$, e.g., $w = 321645 \in S_6$ avoids all fifteen patterns listed in Theorem~\ref{t:pattern-char}.

\begin{theorem}
\label{t:pattern-char}
If the reduced bumpless pipe dream graph $H_w$ for $w$ is a tree, then $w$ avoids the fifteen patterns 1432, 12453, 12534, 13254, 13524, 14253, 14523, 24153, 31524, 35142, 214563, 216345, 354162, 461325, and 465132. Moreover, if $w$ contains any of these fifteen patterns, then $H_w$ contains a $4$-cycle.
\end{theorem}
\begin{proof}
This follows by contrapositive from Propositions~\ref{p:avoid-1432} through~\ref{p:avoid-465132}.
\end{proof}

A natural question after reading Theorem~\ref{t:pattern-char} is whether or not the avoidance of those fifteen patterns is sufficient for the corresponding graph to be a tree. We have seen this to be true for all $w \in S_n$ with $n \leq 9$ using {\sc SageMath}, but we have yet to prove that this remains the case in general. 

\begin{problem}
\label{prob:suff}
Determine sufficient conditions on $w$ in order to guarantee $H_w$ is a tree.
\end{problem}

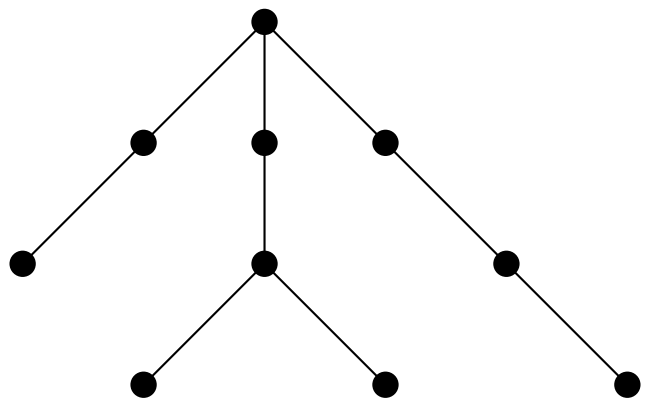
\begin{figure}[h!]
    \centering
\[
\begin{tikzpicture}[scale=0.8]
\coordinate (1) at (0,0) {};
\coordinate (2) at (-1,-1) {};
\coordinate (3) at (0,-1) {};
\coordinate (4) at (1,-1) {};
\coordinate (5) at (-2,-2) {};
\coordinate (6) at (0,-2) {};
\coordinate (7) at (2,-2) {};
\coordinate (8) at (-1,-3) {};
\coordinate (9) at (1,-3) {};
\coordinate (10) at (3,-3) {};

\draw (1) -- (2);
\draw (1) -- (3);
\draw (1) -- (4);
\draw (2) -- (5);
\draw (3) -- (6);
\draw (4) -- (7);
\draw (6) -- (8);
\draw (6) -- (9);
\draw (7) -- (10);

\draw[fill](0,0) circle (0.1);
\draw[fill](-1,-1) circle (0.1);
\draw[fill](0,-1) circle (0.1);
\draw[fill](1,-1) circle (0.1);
\draw[fill](-2,-2) circle (0.1);
\draw[fill](0,-2) circle (0.1);{};
\draw[fill](2,-2) circle (0.1);
\draw[fill](-1,-3) circle (0.1);
\draw[fill](1,-3) circle (0.1);
\draw[fill](3,-3) circle (0.1);

\end{tikzpicture}
\]
\caption{The reduced bumpless pipe dream tree $H_{321645}$.}
\label{fig:tree}
\end{figure}

One might also wonder why we chose to use only minimal droops as edges in $H_w$. We discuss this choice and make connections to other researcher's work in Section~\ref{s:conc}, but we will say just a few words here. The original motivation for this project came from Rubey~\cite{RUBEY}, who defined a partial order on a set of related objects, called {\em reduced pipe dreams}. Reduced pipe dreams were introduced in~\cite{bergeron1993rc} and they are in bijection with reduced bumpless pipe dreams~\cite{gao2021canonical}. Similar to how $\text{BPD}(w)$ can be generated from droop moves, the set of reduced pipe dreams $\text{PD}(w)$ for a permutation $w$ can be generated from what are called {\em chute moves}. Rubey~\cite{RUBEY} introduced the partially ordered set (poset) $C_w = (\text{PD}(w), \leq_{\text{chute}})$, called the {\em chute move poset}, where $P_1 \leq_{\text{chute}} P_2$ if and only if $P_2$ is obtained from $P_1$ via a chute move. Recently, in~\cite{defant} the authors prove Rubey's conjecture: the chute move poset $C_w$ is a lattice for all $w$. Originally, our plan was to introduce the analogous {\em droop move poset} $D_w = (\text{BPD}(w),\leq_{\text{droop}})$. However, droop moves are not transitively closed, e.g., look at Figure~\ref{fig:H1432}. But, as we show in Section~\ref{s:prelim}, every droop move is really a sequence of min-droops, so we should consider $D_w' = (\text{BPD}(w),\leq_{\text{min-droop}})$ instead. \\

The first thing to check is whether or not the Gao--Huang bijection $\varphi$ between $\text{BPD}(w)$ and $\text{PD}(w)$ extends to give a poset isomorphism between $C_w$ and $D_w'$. Alas, it does not, as can be seen by checking the example $w = 1432 \in S_4$. In fact, the poset $D_w'$ need not have a unique maximal element (e.g., $D_{2143}'$ has two maximal elements). Perhaps it would be interesting to characterize which permutations $w$ have that $D_w'$ is a lattice. 

\begin{problem}
\label{prob:poset}
Investigate the min-droop move poset $D_w' = (\text{BPD}(w), \leq_{\text{min-droop}})$, where for any two  $B_1,B_2 \in \text{BPD}(w)$ we have $B_1 \leq_{\text{min-droop}} B_2$ if and only if $B_2$ is obtained from $B_1$ via a min-droop. In particular, characterize the permutations $w$ for which $D_w'$ is a lattice.
\end{problem}

So, we know that $D_w'$ and $C_w$ are very different as posets, which is further testament to just how different that $\text{PD}(w)$ and $\text{BPD}(w)$ really are. Rather than pursue $D_w'$, we decided to ignore the direction on the min-droops and just study the graph $H_w$ instead. \\

\noindent {\bf Organization:} In Section~\ref{s:prelim} we provide the necessary background and definitions, as well as proving some technical lemmas which are the crux of our arguments for the main results in Section~\ref{s:nec-to-avoid}. We conclude with some final remarks and future directions in Section~\ref{s:conc}. \\


\section{Preliminaries}
\label{s:prelim}

For $m,n \in \Z$ with $m \leq n$, we write $[m,n]$ for the set $\{i \in \Z : m \leq i \leq n\}$. We write $[n]$ as shorthand for $[1,n]$. Fix $n$ to be a positive integer throughout the remainder of this section.

\subsection{Permutations}
Let $S_n$ denote the set of bijections $w: [n] \rightarrow [n]$. The elements of $S_n$ are called {\em permutations}. We write a permutation $w \in S_n$ as 
\[
w = w(1) w(2) \ldots w(n),
\] referred to as {\em one-line notation}. For example, $1432 \in S_4$ is the permutation where $w(1) = 1$, $w(2) = 4$, $w(3) = 3$, and $w(4) = 2$. \\

Let $m$ be a positive integer with $m \leq n$. For $w \in S_n$ and $u \in S_m$, we say that $w$ {\em contains the pattern} $u$ if there are indices $1 \leq i_1 \leq \ldots \leq i_m \leq n$ such that the images $w(i_1), \ldots, w(i_m)$ are in the same relative order as $u$. If $w$ does not contain $u$ as a pattern, we say that $w$ {\em avoids} $u$. For example, the permutation 
\[
w = 1 \; \overline{3} \; 5 \; \overline{4} \; 7 \; 6 \; \overline{2} \in S_7
\] 
contains $u = 231 \in S_3$ as a pattern (e.g., the $342$) but avoids $v = 312$. Note that $w$ actually contains eight different instances of the pattern $u$. We call the number of occurrences of a pattern $u$ in a permutation $w$ the {\em packing density of} $u$ {\em in} $w$, and denote this number by $p_u(w)$. Hence $p_{v}(w) = 0$ and one can check that $p_{132}(w) = 11$. \\


\subsection{Bumpless Pipe Dreams}
A \emph{reduced bumpless pipe dream} (BPD) {\em of size} $n$ is a filling of an $n \times n$ grid using tiles of the form
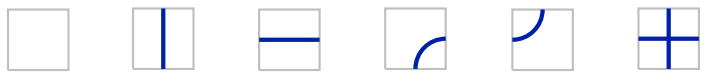
\[
 \nowire \quad \vwire \quad \hwire \quad \are \quad \jay \quad \cross
\]
in such a way that $n$ pipes enter from the bottom and exit to the right, with no pair of pipes crossing each other more than once. We write $\text{BPD}(n)$ for the set of all reduced bumpless pipe dreams of size $n$. We use matrix coordinates to index the $n \times n$ grid of $B \in \text{BPD}(n)$. That is, row indices increase from top to bottom whilst column indices increase from left to right. For $i,j \in [n]$, we use $B(i,j)$ to denote the tile located in position $(i,j)$ of $B$. \\

The pipes are labeled $1$-to-$n$ from left-to-right across the bottom, and the order in which they exit from top-to-bottom determines its {\em associated permutation} $\pi(B) \in S_n$. 
Let $\BPD{w}$ denote the set of all reduced bumpless pipe dreams with associated permutation $w$. 

\subsection{Active Rectangles}
We now introduce our main tool: {\em active rectangles}. We got the name from Weigandt~\cite{weigandt2025}, however, their definition does not allow for $\inlineare$ tiles within the shaded region depicted below in Definition~\ref{def:active}. In this work, we find it convenient to allow for $\inlineare$ tiles inside the rectangles we consider. Note that every active rectangle (by our definition) contains a sub-rectangle that is active by Weigandt's definition. 

\begin{definition}
\label{def:active}
    A rectangle $[a,b] \times [c,d]$ of tiles in a bumpless pipe dream $B$ is called \emph{active} if 
    \begin{minipage}{0.49 \textwidth}







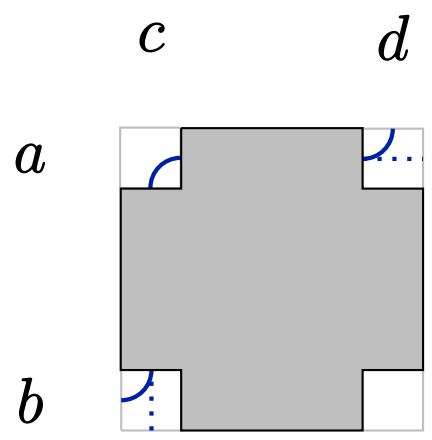
    \begin{itemize}
        \item[(i)] $B(a,c) = \inlineare$,
        \item[(ii)] $B(b,d) = \inlinenowire$,
        \item[(iii)] $B(a,d) \in \{\inlinejay \; , \; \inlinehwire\}$, 
        \item[(iv)] $B(b,c) \in \{\inlinejay \; , \; \inlinevwire\}$, and
        \item[(v)] $B(i,j) \neq \inlinejay$ for all other $i,j \in [n]$.
    \end{itemize}
    \end{minipage}
    \begin{minipage}{0.49 \textwidth}
	\[
	\begin{tikzpicture}[scale=0.4]
		\node at (0,0) {\are};
		\node at (4,-4) {\nowire};
		\node at (0,-4) {\vewire};
		\node at (4,0) {\hewire};
		
		\node at (-2,0) {$a$};
		\node at (-2,-4) {$b$};
		\node at (0,2) {$c$};
		\node at (4,2) {$d$};

		\draw[fill=lightgray] (0.5,0.5) -- (0.5,-0.5) -- (-0.5,-0.5) -- (-0.5,-3.5) -- (0.5,-3.5) -- (0.5,-4.5) -- (3.5,-4.5) -- (3.5,-3.5) -- (4.5,-3.5) -- (4.5,-0.5) -- (3.5,-0.5) -- (3.5,0.5) -- (0.5,0.5);
	\end{tikzpicture}
	\]
	\end{minipage}
\end{definition}

The following two results follow readily from definitions and are used frequently in the rest of this article.

\begin{proposition}
\label{prop:border}
If $R = [a,b] \times [c,d]$ is an active rectangle on a reduced bumpless pipe dream $B$, then the tiles on the border of $R$ are constrained. In particular, we have:


\begin{itemize}
\item[(a)] $B(a,j) \in \{\inlinehwire \; , \; \inlinecross\}$ for all $c < j < d$, \hfill (top border)
\item[(b)] $B(i,c) \in \{\inlinevwire \; , \; \inlinecross\}$ for all $a < i < b$, \hfill (left border)
\item[(c)] $B(b,j) \in \{\inlinenowire \; , \; \inlinevwire\}$ for all $c < j < d$, and \hfill (bottom border)
\item[(d)] $B(i,d) \in \{\inlinenowire \; , \; \inlinehwire\}$ for all $a < i < b$. \hfill (right border)

\end{itemize}
\end{proposition}
\begin{proposition}
\label{prop:sub}
    If $R = [a,b] \times [c,d]$ is an active rectangle on a reduced bumpless pipe dream $B$ and $B(a',c') = \inlineare$ for some $(a',c')$ with $a < a' < b$ and $c < c' < d$, then the subset
    \[
    R' = [a',b] \times [c',d]
    \]
    is also an active rectangle on $B$.
\end{proposition}

\subsection{Droop Moves} 
The set $\text{BPD}(w)$ is connected via local moves, called ``droop moves''. A {\em droop move} on some $B \in \text{BPD}(w)$ is a reassignment of the tiles matching one of the following forms:
    
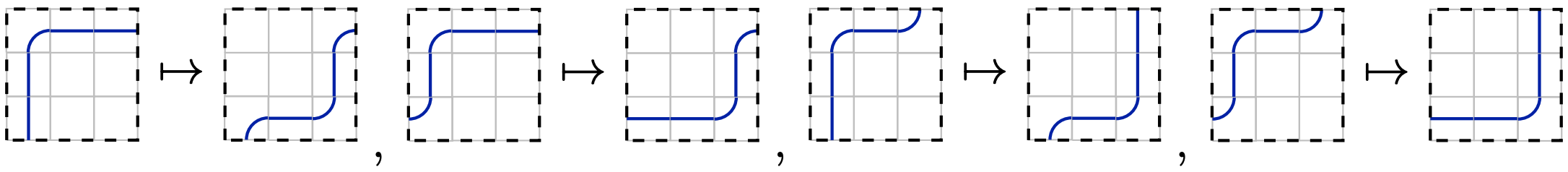
\[
\begin{tikzpicture}[scale= .4]
\node at (0,0) {\vwire};
\node at (1,0) {\nowire};
\node at (2,0) {\nowire};

\node at (0,1) {\vwire};
\node at (1,1) {\nowire};
\node at (2,1) {\nowire};

\node at (0,2) {\are};
\node at (1,2) {\hwire};
\node at (2,2) {\hwire};

\draw[thick, dashed] (-0.5,-0.5) -- (-0.5,2.5) -- (2.5,2.5) -- (2.5,-0.5) -- (-0.5,-0.5);

\node at (5,0) {\are};
\node at (6,0) {\hwire};
\node at (7,0) {\jay};

\node at (5,1) {\nowire};
\node at (6,1) {\nowire};
\node at (7,1) {\vwire};

\node at (5,2) {\nowire};
\node at (6,2) {\nowire};
\node at (7,2) {\are};

\node at (3.5, 1) {$\mapsto$};

\draw[thick, dashed] (4.5,-0.5) -- (4.5,2.5) -- (7.5,2.5) -- (7.5,-0.5) -- (4.5,-0.5);

\end{tikzpicture},
\begin{tikzpicture}[scale= .4]
\node at (0,0) {\jay};
\node at (1,0) {\nowire};
\node at (2,0) {\nowire};

\node at (0,1) {\vwire};
\node at (1,1) {\nowire};
\node at (2,1) {\nowire};

\node at (0,2) {\are};
\node at (1,2) {\hwire};
\node at (2,2) {\hwire};

\draw[thick, dashed] (-0.5,-0.5) -- (-0.5,2.5) -- (2.5,2.5) -- (2.5,-0.5) -- (-0.5,-0.5);

\node at (5,0) {\hwire};
\node at (6,0) {\hwire};
\node at (7,0) {\jay};

\node at (5,1) {\nowire};
\node at (6,1) {\nowire};
\node at (7,1) {\vwire};

\node at (5,2) {\nowire};
\node at (6,2) {\nowire};
\node at (7,2) {\are};

\node at (3.5, 1) {$\mapsto$};

\draw[thick, dashed] (4.5,-0.5) -- (4.5,2.5) -- (7.5,2.5) -- (7.5,-0.5) -- (4.5,-0.5);
\end{tikzpicture},
\begin{tikzpicture}[scale= .4]
\node at (0,0) {\vwire};
\node at (1,0) {\nowire};
\node at (2,0) {\nowire};

\node at (0,1) {\vwire};
\node at (1,1) {\nowire};
\node at (2,1) {\nowire};

\node at (0,2) {\are};
\node at (1,2) {\hwire};
\node at (2,2) {\jay};

\draw[thick, dashed] (-0.5,-0.5) -- (-0.5,2.5) -- (2.5,2.5) -- (2.5,-0.5) -- (-0.5,-0.5);

\node at (5,0) {\are};
\node at (6,0) {\hwire};
\node at (7,0) {\jay};

\node at (5,1) {\nowire};
\node at (6,1) {\nowire};
\node at (7,1) {\vwire};

\node at (5,2) {\nowire};
\node at (6,2) {\nowire};
\node at (7,2) {\vwire};

\node at (3.5, 1) {$\mapsto$};

\draw[thick, dashed] (4.5,-0.5) -- (4.5,2.5) -- (7.5,2.5) -- (7.5,-0.5) -- (4.5,-0.5);

\end{tikzpicture},
\begin{tikzpicture}[scale= .4]
\node at (0,0) {\jay};
\node at (1,0) {\nowire};
\node at (2,0) {\nowire};

\node at (0,1) {\vwire};
\node at (1,1) {\nowire};
\node at (2,1) {\nowire};

\node at (0,2) {\are};
\node at (1,2) {\hwire};
\node at (2,2) {\jay};

\draw[thick, dashed] (-0.5,-0.5) -- (-0.5,2.5) -- (2.5,2.5) -- (2.5,-0.5) -- (-0.5,-0.5);

\node at (5,0) {\hwire};
\node at (6,0) {\hwire};
\node at (7,0) {\jay};

\node at (5,1) {\nowire};
\node at (6,1) {\nowire};
\node at (7,1) {\vwire};

\node at (5,2) {\nowire};
\node at (6,2) {\nowire};
\node at (7,2) {\vwire};

\node at (3.5, 1) {$\mapsto$};

\draw[thick, dashed] (4.5,-0.5) -- (4.5,2.5) -- (7.5,2.5) -- (7.5,-0.5) -- (4.5,-0.5);

\end{tikzpicture}
\]
where unpictured pipes may pass through the region, but no $\inlineare$ or $\inlinejay$ tiles appear within. We can use active rectangles $R = [a,b] \times [c,d]$ to describe those regions depicted above in droop moves, where $a < b$ give the row indices under consideration and $c < d$ give the column indices. Hence, the top left corner is a $\inlineare$ tile in location $(a,c)$ and the bottom right corner is a $\inlinenowire$ tile in location $(b,d)$, and there are no $\inlineare$ or $\inlinejay$ tiles in $R - \{(a,c),(a,d),(b,c)\}$. In this case, we write $(a,c) \searrow (b,d)$ to denote the droop move. The following result follows readily from definitions.

\begin{proposition}
\label{prop:min}
    If $R = [a,b] \times [c,d]$ is a rectangle on a reduced bumpless pipe dream $B$, then the droop move $(a,c) \searrow (b,d)$ is valid if and only if $R$ is active and $B(i,j) \neq \inlineare$ for all $(i,j) \neq (a,c)$.
\end{proposition}

So, we can think of droop moves as taking place inside of active rectangles that do not contain smaller active rectangles inside them. 

\subsection{Minimal Droops}
A droop move $(a,b) \searrow (c,d)$ is called {\em minimal} (or a {\em min-droop} for short) if there are no other $\inlinenowire$ tiles in the active rectangle $[a,b] \times [c,d]$. We now show that droop moves can be decomposed into a sequence of min-droops.

\begin{lemma}
\label{l:droop-imply-min-droop}
Every droop move is a sequence of min-droop moves. 
\end{lemma}
\begin{proof}
	Let $R = [a,b] \times [c,d]$ be a rectangle of tiles on a reduced bumpless pipe dream and suppose $(a,c) \searrow (b,d)$ is a valid droop move. By Proposition~\ref{prop:min}, this implies $R$ is active and there are no $\inlineare$ or $\inlinejay$ tiles in the region $R' = R - \{(a,c),(a,d), (b,c)\}$. Thus, all tiles in $R'$ come from the set $\{\inlinenowire \; , \inlinevwire \; , \inlinehwire \; , \inlinecross \}$. We proceed by induction on the number of $\inlinenowire$ tiles in $R-\{(b,d)\}$. For our base case, suppose there are no $\inlinenowire$ tiles in $R-\{(b,d)\}$. Then, by definition, $(a,c) \searrow (b,d)$ is a min-droop move as desired. Now suppose there are $n > 0$ total $\inlinenowire$ tiles in $R - \{(b,d)\}$, and that the statement holds whenever there are fewer than $n$ total $\inlinenowire$ tiles in a given droop region. Without loss of generality, pick any such $\inlinenowire$ tile and suppose it is in location $(i,j)$, with $(i,j) \neq (b,d)$. 

    \underline{Case 1:} Suppose that $i < b$. Observe that, since $R'$ has no $\inlineare$ or $\inlinejay$ tiles, by Proposition~\ref{prop:border}(d) we have that $B(i,d) = \inlinenowire$. Therefore, the rectangle $R_1 = [a,i] \times [c,d]$ is active and $R_1 - \{(a,c)\}$ contains no $\inlineare$ tiles, so then $(a,c) \searrow (i,d)$ is a valid droop move by Proposition~\ref{prop:min}. Since $R_1$ contains fewer than $n$ total $\inlinenowire$ tiles, then by induction, the droop move $(a,c) \searrow (i,d)$ decomposes into a sequence of min-droops that result in a reduced bumpless pipe dream $B'$ with $B'(i,c) = \inlineare$ and $B'(i,d) = \inlinejay$, and all other tiles outside of $R_1$ agree on both $B$ and $B'$. Now consider the rectangle $R_2 = [i,b] \times [c,d]$ on $B'$. Observe that $R_2$ is active, there are no $\inlineare$ tiles in $R_2 - \{(i,c)\}$, and there are fewer than $n$ total $\inlinenowire$ tiles on $R_2$. So, the droop move $(i,c) \searrow (b,d)$ is valid on $B'$ and, again by induction, we obtain a sequence of min-droops, starting at $B'$ and resulting in a reduced bumpless pipe dream $B''$. Concatenating these two sequences together (from $B$ to $B'$ and then $B'$ to $B''$) gives the desired decomposition of $(a,c) \searrow (b,d)$ into min-droops. 

    \underline{Case 2:} Suppose that $i = b$. Then we must have that $j < d$ and now we appeal to Proposition~\ref{prop:border}(c) to deduce that $B(b,j) = \inlinenowire$. Then let $R_1 = [a,b] \times [c,j]$ and follow the analogous steps as we did in Case 1. In short, the droop move $(a,c) \searrow (b,j)$ is valid, which gives a sequence of min-droops to a reduced bumpless pipe dream $B'$, where $B'(a,j) = \inlineare$ and $B'(b,j) = \inlinejay$. Now one can verify that $R_2 = [a,b] \times [j,d]$ is an active rectangle, the droop move $(a,j) \searrow (b,d)$ is valid on $B'$, and $R_2$ has fewer than $n$ total $\inlinenowire$ tiles. Hence, after applying the inductive hypothesis to $R_2$ on $B'$ and concatenating sequences, we again obtain the desired decomposition of $(a,c) \searrow (b,d)$ into min-droops. 
 \end{proof}

\subsection{Graph Structure}

\begin{definition}
For a permutation $w$, the {\em reduced bumpless pipe dream graph} for $w$ is the graph given by $H_w = (V_w,E_w)$ where $V_w = \BPD{w}$ and $\{B_1,B_2\} \in E_w$ if and only if $B_1$ and $B_2$ are related by a min-droop.
\end{definition}

\begin{corollary}
    The reduced bumpless pipe dream graph $H_w$ is connected for any $w$.
\end{corollary}

Examples of the reduced bumpless pipe dream graphs $H_{1432}$ and $H_{321645}$ are shown in Figures~\ref{fig:H1432} and~\ref{fig:tree}, respectively. We now establish a technical lemma that is used to shorten our arguments for the main results in the next section. 
\begin{lemma}
\label{l:min-droop-rectangle}
If a rectangle $R = [a,b] \times [c,d]$ in a bumpless pipe dream $B$ is active, then there is a sequence of min-droop moves $B = B_0 \to B_1 \to \ldots \to B_\ell$, with $\ell \geq 1$, where each distinct $B_i$ and $B_j$ differ only by tiles in $R$ and $B_\ell(a,c) = \inlinenowire$.
\end{lemma}
\begin{proof}
	We proceed by using induction on the number of $\inlineare$ tiles in the region $R' = R - \{(a,c)\}$. For our base case, suppose there are no $\inlineare$ tiles in $R$. Since $R$ is active, it follows that the droop move $(a,c) \searrow (b,d)$ is available. This droop move only modifies tiles in $R$, so by Lemma~\ref{l:droop-imply-min-droop}, we get a sequence of min-droop moves as advertised.
	
	Now suppose there are $n > 0$ total $\inlineare$ tiles in $R$, and assume that the statement holds whenever there are fewer than $n$ total $\inlineare$ tiles in a given region. Pick any such $\inlineare$ tile and suppose it is in cell $(i,j)$ for some $a < i < b$ and $c < j < d$. (Note that these inequalities must be strict by Proposition~\ref{prop:border}.) Consider the rectangle $R_1 = [i,b] \times [j,d]$. Since this region is active (by Proposition~\ref{prop:sub}) and has at most $n-1$ total $\inlineare$ tiles, by induction we know there is a sequence of min-droop moves that results in a reduced bumpless pipe dream $B_\ell$ with $B_\ell(i,j) = \inlinenowire$. Now consider the rectangle $R_2 = [a,i] \times [c,j]$ on $B_\ell$. Observe that $R_2$ is active. Moreover, since $R_2$ has at most $n-1$ total $\inlineare$ tiles, by induction we get a sequence of min-droop moves that begins at $B_\ell$ and terminates at some $B_\ell'$ with $B_\ell'(a,c) = \inlinenowire$. Concatenating these sequences gives the result.
\end{proof}
The following two corollaries are immediate from Lemma~\ref{l:min-droop-rectangle}.
\begin{corollary}
\label{cor:min-droop}
Every active rectangle contains a min-droop.
\end{corollary}
\begin{corollary}
\label{cor:main}
If any reduced bumpless pipe dream for $w$ contains two active rectangles that are disjoint or intersect only at the bottom-left corner of one rectangle and the top-right corner of the other, then $H_w$ contains a $4$-cycle.  
\end{corollary}
\begin{proof}
Let $w$ be a permutation and suppose $B \in \text{BPD}(w)$ has two active rectangles $R_1$ and $R_2$ that are either disjoint, or intersect each other at the bottom-left corner of $R_1$ and top-right corner of $R_2$. Since they are both active, this forces that tile to be a $\inlinejay$ tile, if it exists. By Corollary~\ref{cor:min-droop}, the rectangle $R_1$ contains a min-droop. Let $B_1$ be the reduced bumpless pipe dream for $w$ obtained by performing this min-droop on $B$. Consider $R_2$ on $B_1$. Every tile in $R_2$ is the same on $B_1$ as it is on $B$, except for possibly the top-right corner. If $R_1$ and $R_2$ are disjoint then this tile is unchanged. If $R_1$ and $R_2$ intersect, then on $B$ the tile must have been a $\inlinejay$ tile as $R_1$ and $R_2$ are both active. Hence, on $B_1$, this tile is now either a $\inlinejay$ still, or it is a $\inlinehwire$ now. Regardless, $R_2$ remains active on $B_1$. Therefore, by Corollary~\ref{cor:min-droop}, the rectangle $R_2$ contains a min-droop on $B_1$. Let $B_3$ be the result of performing this min-droop on $B_1$. Observe that we could instead have chosen to perform the min-droop in $R_2$ on $B$ first, to obtain some $B_2 \in \text{BPD}(w)$, and then have that $R_1$ remains active on $B_2$ and perform its min-droop on $B_2$ to get $B_3$. Thus $\{B,B_1,B_2,B_3\} \subseteq V(H_w)$ form a $4$-cycle in $H_w$.
\end{proof}

\section{Main Results}
\label{s:nec-to-avoid}
In this section we establish the necessity of $w$ to avoid each of the fifteen patterns whenever $H_w$ is a tree. We do this for each pattern individually. The proof strategy for each pattern is the same, so we will give the general outline here. We start by supposing $w$ contains a certain pattern, then we work to identify a reduced bumpless pipe dream for $w$ which contains two active rectangles that are either disjoint, or meet at the allowed corners. From there, we can appeal to Corollary~\ref{cor:main} to produce a $4$-cycle in $H_w$, thus proving the claimed necessity via contrapositive. \\

Throughout this section, we call $\inlineare$ and $\inlinejay$ tiles {\em elbows}.
\begin{proposition}
\label{p:avoid-1432}
	If $w$ contains a 1432 pattern, then $H_w$ contains a $4$-cycle.
\end{proposition}
\begin{proof}
Let $k_1 < k_2 < k_3 < k_4$ be an occurrence of a 1432 pattern in $w$, so that 
	\[\
	w(k_1) < w(k_4) < w(k_3) < w(k_2).
	\]
Consider $\overline{B}_w$, the Rothe BPD for $w$. Use $j_1 < j_2 < j_3 < j_4$ to label the four rows of $\overline{B}_w$ corresponding to the unique $\inlineare$ tiles of this pattern. Observe that $k_1 < k_2 < k_3 < k_4$ label the four columns of $\overline{B}_w$ corresponding to the unique $\inlineare$ tiles of this pattern. Then we have 
\[
\overline{B}_w(j_2,k_2) = \inlinenowire \; , \quad \overline{B}_w(j_3,k_2) = \inlinenowire \; , \quad \text{and} \quad
\overline{B}_w(j_2,k_3) = \inlinenowire
\]
due to the inversions $(k_2,k_4)$, $(k_2,k_3)$, and $(k_3,k_4)$, respectively, in $w$. Consider the rectangle
\[
R = [j_1,j_2] \times [k_1,k_2] \subseteq \overline{B}_w.
\]
In $R$, the top-left corner is $\overline{B}_w(j_1,k_1) = \inlineare$ and the bottom-right corner is $\overline{B}_w(j_2,k_2) = \inlinenowire$. Moreover, there are no $\inlinejay$ tiles in $R$ since it is a rectangle on $\overline{B}_w$. If there is a $\inlineare$ tile in $R - \{(j_1,k_1)\}$, it would be in some row $j_1'$ with $j_1 < j_1' < j_2$ and some column $k_1'$ with $k_1 < k_1' < k_2$. In this case, we could select $k_1' \; k_2 \; k_3 \; k_4$ as our occurrence of a 1432 pattern instead. Hence, it suffices to assume there are no elbows in $R$, which implies that the droop move $(j_1,k_1) \searrow (j_2,k_2)$ is valid. Let $B$ denote the BPD for $w$ obtained by performing this droop move on $\overline{B}_w$. Note that the tiles in $B$ and $\overline{B}_w$ only differ inside of $R$, and more specifically, we have the following reassignment of tiles on the corners of $R$:
\[
B(j_1,k_2) = \inlineare \; , \quad B(j_2,k_1) = \inlineare \; , \quad \text{and} \quad B(j_2,k_2) = \inlinejay.
\]
Consider the following two rectangles on $B$:
\[
R_1 = [j_1,j_2] \times [k_2,k_3] \subseteq B \quad \text{and} \quad R_2 = [j_2,j_3] \times [k_1,k_2] \subseteq B.
\]
Observe that $R_1 \cap R_2 = \{(j_2,k_2)\}$. Moreover, we claim that both $R_1$ and $R_2$ are active. Indeed, the only $\inlinejay$ tile in $B$ is $B(j_2,k_2)$, which is the bottom-left corner of $R_1$ and top-right corner of $R_2$. For $R_1$, the top-left corner is $B(j_1,k_2) = \inlineare$, the top-right corner of $R_1$ is $B(j_1,k_3)= \inlinehwire$, and the bottom-right corner is $B(j_2,k_3)=\inlinenowire$. Hence, $R_1$ is active. Likewise, for $R_2$, the top-left corner is $B(j_2,k_1) = \inlineare$, the bottom-left corner is $B(j_3,k_1) = \inlinevwire$, and the bottom-right corner is $B(j_3,k_2) = \inlinenowire$. Hence, $R_2$ is active. Therefore, by Corollary~\ref{cor:main}, we have that $H_w$ contains a 4-cycle as advertised.
\end{proof}
\begin{proposition}
\label{p:avoid-12453}
	If $w$ contains a 12453 pattern, then $H_w$ contains a $4$-cycle.
\end{proposition}
\begin{proof}
Let $k_1 < \ldots < k_5$ be an occurrence of a 12453 pattern in $w$, so that
\[
w(k_1) < w(k_2) < w(k_4) < w(k_5) < w(k_3).
\]
Consider $\overline{B}_w$, the Rothe BPD for $w$. Use $j_1 < \ldots < j_5$ to label the five rows of $\overline{B}_w$ corresponding to the unique $\inlineare$ tiles of this pattern. Observe that $k_1 < \ldots < k_5$ label the five columns of $\overline{B}_w$ corresponding to the unique $\inlineare$ tiles of this pattern. Then we have 
\[
\overline{B}_w(j_3,k_3) = \inlinenowire \quad \text{and} \quad
\overline{B}_w(j_4,k_3) = \inlinenowire 
\]
due to the inversions $(k_3,k_5)$ and $(k_4,k_5)$, respectively, in $w$. Consider the rectangle
\[
R = [j_2,j_3] \times [k_2,k_3] \subseteq \overline{B}_w.
\]
In $R$, the top-left corner is $\overline{B}_w(j_2,k_2) = \inlineare$ and the bottom-right corner is $\overline{B}_w(j_3,k_3) = \inlinenowire$. Moreover, there are no $\inlinejay$ tiles in $R$ since it is a rectangle on $\overline{B}_w$. If there is a $\inlineare$ tile in $R - \{(j_2,k_2)\}$, it would be in some row $j_2'$ with $j_2 < j_2' < j_3$ and some column $k_2'$ with $k_2 < k_2' < k_3$. In this case, we could select $k_2 \; k_2' \; k_3 \; k_4 \; k_5$ as our occurrence of a 12453 pattern instead. Hence, it suffices to assume there are no elbows in $R$, which implies that the droop move $(j_2,k_2) \searrow (j_3,k_3)$ is valid. Let $B$ denote the BPD for $w$ obtained by performing this droop move on $\overline{B}_w$. Note that the tiles in $B$ and $\overline{B}_w$ only differ inside of $R$, and more specifically, we have the following reassignment of tiles on the corners of $R$:
\[
B(j_3,k_2) = \inlineare \; , \quad B(j_2,k_3) = \inlineare \; , \quad \text{and} \quad B(j_3,k_3) = \inlinejay.
\]
Consider the following two rectangles on $B$:
\[
R_1 = [j_1,j_2] \times [k_1,k_2] \subseteq B \quad \text{and} \quad R_2 = [j_3,j_4] \times [k_2,k_3] \subseteq B.
\]
Observe that $R_1$ and $R_2$ are disjoint since $j_2 < j_3$. Moreover, we claim that both $R_1$ and $R_2$ are active. Indeed, the only $\inlinejay$ tile in $B$ is $B(j_3,k_3)$, which is the top-right corner of $R_2$. For $R_1$, the top-left corner is $B(j_1,k_1) = \inlineare$, the top-right corner of $R_1$ is $B(j_1,k_2)= \inlinehwire$, and the bottom-right corner is $B(j_2,k_2)=\inlinenowire$. Hence, $R_1$ is active. Likewise, for $R_2$, the top-left corner is $B(j_3,k_2) = \inlineare$, the bottom-left corner is $B(j_4,k_2) = \inlinevwire$, and the bottom-right corner is $B(j_4,k_3) = \inlinenowire$. Hence, $R_2$ is active. Therefore, by Corollary~\ref{cor:main}, we have that $H_w$ contains a 4-cycle as advertised.
\end{proof}
\begin{proposition}
\label{p:avoid-12534}
	If $w$ contains a 12534 pattern, then $H_w$ contains a $4$-cycle.
\end{proposition}
\begin{proof}
Let $k_1 < \ldots < k_5$ be an occurrence of a 12534 pattern in $w$, so that
\[
w(k_1) < w(k_2) < w(k_5) < w(k_3) < w(k_4).
\]
Consider $\overline{B}_w$, the Rothe BPD for $w$. Use $j_1 < \ldots < j_5$ to label the five rows of $\overline{B}_w$ corresponding to the unique $\inlineare$ tiles of this pattern. Observe that $k_1 < \ldots < k_5$ label the five columns of $\overline{B}_w$ corresponding to the unique $\inlineare$ tiles of this pattern. Then we have 
\[
\overline{B}_w(j_3,k_3) = \inlinenowire \quad \text{and} \quad
\overline{B}_w(j_3,k_4) = \inlinenowire 
\]
due to the inversions $(k_3,k_4)$ and $(k_3,k_5)$, respectively, in $w$. Consider the rectangle
\[
R = [j_2,j_3] \times [k_2,k_3] \subseteq \overline{B}_w.
\]
In $R$, the top-left corner is $\overline{B}_w(j_2,k_2) = \inlineare$ and the bottom-right corner is $\overline{B}_w(j_3,k_3) = \inlinenowire$. Moreover, there are no $\inlinejay$ tiles in $R$ since it is a rectangle on $\overline{B}_w$. If there is a $\inlineare$ tile in $R - \{(j_2,k_2)\}$, it would be in some row $j_2'$ with $j_2 < j_2' < j_3$ and some column $k_2'$ with $k_2 < k_2' < k_3$. In this case, we could select $k_2 \; k_2' \; k_3 \; k_4 \; k_5$ as our occurrence of a 12534 pattern instead. Hence, it suffices to assume there are no elbows in $R$, which implies that the droop move $(j_2,k_2) \searrow (j_3,k_3)$ is valid. Let $B$ denote the BPD for $w$ obtained by performing this droop move on $\overline{B}_w$. Note that the tiles in $B$ and $\overline{B}_w$ only differ inside of $R$, and more specifically, we have the following reassignment of tiles on the corners of $R$:
\[
B(j_3,k_2) = \inlineare \; , \quad B(j_2,k_3) = \inlineare \; , \quad \text{and} \quad B(j_3,k_3) = \inlinejay.
\]
Consider the following two rectangles on $B$:
\[
R_1 = [j_1,j_2] \times [k_1,k_2] \subseteq B \quad \text{and} \quad R_2 = [j_2,j_3] \times [k_3,k_4] \subseteq B.
\]
Observe that $R_1$ and $R_2$ are disjoint since $k_2 < k_3$. Moreover, we claim that both $R_1$ and $R_2$ are active. Indeed, the only $\inlinejay$ tile in $B$ is $B(j_3,k_3)$, which is the bottom-left corner of $R_2$. For $R_1$, the top-left corner is $B(j_1,k_1) = \inlineare$, the top-right corner of $R_1$ is $B(j_1,k_2)= \inlinehwire$, and the bottom-right corner is $B(j_2,k_2)=\inlinenowire$. Hence, $R_1$ is active. Likewise, for $R_2$, the top-left corner is $B(j_2,k_3) = \inlineare$, the top-right corner is $B(j_2,k_4) = \inlinehwire$, and the bottom-right corner is $B(j_3,k_4) = \inlinenowire$. Hence, $R_2$ is active. Therefore, by Corollary~\ref{cor:main}, we have that $H_w$ contains a 4-cycle as advertised.
\end{proof}
\begin{proposition}
\label{p:avoid-13254}
	If $w$ contains a 13254 pattern, then $H_w$ contains a $4$-cycle.
\end{proposition}
\begin{proof}
Let $k_1 < \ldots < k_5$ be an occurrence of a 13254 pattern in $w$, so that
\[
w(k_1) < w(k_3) < w(k_2) < w(k_5) < w(k_4).
\]
Consider $\overline{B}_w$, the Rothe BPD for $w$. Use $j_1 < \ldots < j_5$ to label the five rows of $\overline{B}_w$ corresponding to the unique $\inlineare$ tiles of this pattern. Observe that $k_1 < \ldots < k_5$ label the five columns of $\overline{B}_w$ corresponding to the unique $\inlineare$ tiles of this pattern. Consider the following two rectangles on $\overline{B}_w$:
\[
R_1 = [j_1,j_2] \times [k_1,k_2] \subseteq \overline{B}_w \quad \text{and} \quad R_2 = [j_3,j_4] \times [k_2,k_4] \subseteq \overline{B}_w.
\]
Observe that $R_1$ and $R_2$ are disjoint since $j_2 < j_3$. Moreover, we claim that both $R_1$ and $R_2$ are active. Indeed, there are no $\inlinejay$ tiles in either $R_1$ or $R_2$, since they are rectangles on $\overline{B}_w$. For $R_1$, the top-left corner is $\overline{B}_w(j_1,k_1) = \inlineare$, the top-right corner of $R_1$ is $\overline{B}_w(j_1,k_2)= \inlinehwire$, the bottom-left corner of $R_1$ is $\overline{B}_w(j_2,k_1) = \inlinevwire$, and the bottom-right corner is $\overline{B}_w(j_2,k_2)=\inlinenowire$. Hence, $R_1$ is active. Likewise, for $R_2$, the top-left corner is $\overline{B}_w(j_3,k_2) = \inlineare$, the top-right corner is $\overline{B}_w(j_3,k_4) = \inlinehwire$, the bottom-left corner is $\overline{B}_w(j_4,k_2) = \inlinevwire$, and the bottom-right corner is $B(j_4,k_4) = \inlinenowire$. Hence, $R_2$ is active. Therefore, by Corollary~\ref{cor:main}, we have that $H_w$ contains a 4-cycle as advertised.
\end{proof}

\begin{proposition}
\label{p:avoid-13524}
	If $w$ contains a 13524 pattern, then $H_w$ contains a $4$-cycle.
\end{proposition}
\begin{proof}
Let $k_1 < \ldots < k_5$ be an occurrence of a 13524 pattern in $w$, so that
\[
w(k_1) < w(k_3) < w(k_5) < w(k_2) < w(k_4).
\]
Consider $\overline{B}_w$, the Rothe BPD for $w$. Use $j_1 < \ldots < j_5$ to label the five rows of $\overline{B}_w$ corresponding to the unique $\inlineare$ tiles of this pattern. Observe that $k_1 < \ldots < k_5$ label the five columns of $\overline{B}_w$ corresponding to the unique $\inlineare$ tiles of this pattern. Consider the following two rectangles on $\overline{B}_w$:
\[
R_1 = [j_1,j_2] \times [k_1,k_2] \subseteq \overline{B}_w \quad \text{and} \quad R_2 = [j_2,j_3] \times [k_3,k_4] \subseteq \overline{B}_w.
\]
Observe that $R_1$ and $R_2$ are disjoint since $k_2 < k_3$. Moreover, we claim that both $R_1$ and $R_2$ are active. Indeed, there are no $\inlinejay$ tiles in either $R_1$ or $R_2$, since they are rectangles on $\overline{B}_w$. For $R_1$, the top-left corner is $\overline{B}_w(j_1,k_1) = \inlineare$, the top-right corner of $R_1$ is $\overline{B}_w(j_1,k_2)= \inlinehwire$, the bottom-left corner of $R_1$ is $\overline{B}_w(j_2,k_1) = \inlinevwire$, and the bottom-right corner is $\overline{B}_w(j_2,k_2)=\inlinenowire$. Hence, $R_1$ is active. Likewise, for $R_2$, the top-left corner is $\overline{B}_w(j_3,k_3) = \inlineare$, the top-right corner is $\overline{B}_w(j_3,k_4) = \inlinehwire$, the bottom-left corner is $\overline{B}_w(j_3,k_3) = \inlinevwire$, and the bottom-right corner is $B(j_3,k_4) = \inlinenowire$. Hence, $R_2$ is active. Therefore, by Corollary~\ref{cor:main}, we have that $H_w$ contains a 4-cycle as advertised.
\end{proof}
\begin{proposition}
\label{p:avoid-14253}
	If $w$ contains a 14253 pattern, then $H_w$ contains a $4$-cycle.
\end{proposition}
\begin{proof}
Let $k_1 < \ldots < k_5$ be an occurrence of a 14253 pattern in $w$, so that
\[
w(k_1) < w(k_4) < w(k_2) < w(k_5) < w(k_3).
\]
Consider $\overline{B}_w$, the Rothe BPD for $w$. Use $j_1 < \ldots < j_5$ to label the five rows of $\overline{B}_w$ corresponding to the unique $\inlineare$ tiles of this pattern. Observe that $k_1 < \ldots < k_5$ label the five columns of $\overline{B}_w$ corresponding to the unique $\inlineare$ tiles of this pattern. Consider the following two rectangles on $\overline{B}_w$:
\[
R_1 = [j_1,j_2] \times [k_1,k_2] \subseteq \overline{B}_w \quad \text{and} \quad R_2 = [j_3,j_4] \times [k_2,k_3] \subseteq \overline{B}_w.
\]
Observe that $R_1$ and $R_2$ are disjoint since $j_2 < j_3$. Moreover, we claim that both $R_1$ and $R_2$ are active. Indeed, there are no $\inlinejay$ tiles in either $R_1$ or $R_2$, since they are rectangles on $\overline{B}_w$. For $R_1$, the top-left corner is $\overline{B}_w(j_1,k_1) = \inlineare$, the top-right corner of $R_1$ is $\overline{B}_w(j_1,k_2)= \inlinehwire$, the bottom-left corner of $R_1$ is $\overline{B}_w(j_2,k_1) = \inlinevwire$, and the bottom-right corner is $\overline{B}_w(j_2,k_2)=\inlinenowire$. Hence, $R_1$ is active. Likewise, for $R_2$, the top-left corner is $\overline{B}_w(j_3,k_2) = \inlineare$, the top-right corner is $\overline{B}_w(j_3,k_3) = \inlinehwire$, the bottom-left corner is $\overline{B}_w(j_4,k_2) = \inlinevwire$, and the bottom-right corner is $B(j_4,k_3) = \inlinenowire$. Hence, $R_2$ is active. Therefore, by Corollary~\ref{cor:main}, we have that $H_w$ contains a 4-cycle as advertised.
\end{proof}

\begin{proposition}
\label{p:avoid-14523}
	If $w$ contains a 14523 pattern, then $H_w$ contains a $4$-cycle.
\end{proposition}
\begin{proof}
Let $k_1 < \ldots < k_5$ be an occurrence of a 14523 pattern in $w$, so that
\[
w(k_1) < w(k_4) < w(k_5) < w(k_2) < w(k_3).
\]
Consider $\overline{B}_w$, the Rothe BPD for $w$. Use $j_1 < \ldots < j_5$ to label the five rows of $\overline{B}_w$ corresponding to the unique $\inlineare$ tiles of this pattern. Observe that $k_1 < \ldots < k_5$ label the five columns of $\overline{B}_w$ corresponding to the unique $\inlineare$ tiles of this pattern. Consider the rectangle
\[
R = [j_1,j_2] \times [k_1,k_2] \subseteq \overline{B}_w.
\]
In $R$, the top-left corner is $\overline{B}_w(j_1,k_1) = \inlineare$ and the bottom-right corner is $\overline{B}_w(j_2,k_2) = \inlinenowire$. Moreover, there are no $\inlinejay$ tiles in $R$ since it is a rectangle on $\overline{B}_w$. If there is a $\inlineare$ tile in $R - \{(j_1,k_1)\}$, it would be in some row $j_1'$ with $j_1 < j_1' < j_2$ and some column $k_1'$ with $k_1 < k_1' < k_2$. In this case, we could select $k_1' \; k_2 \; k_3 \; k_4 \; k_5$ as our occurrence of a 14523 pattern instead. Hence, it suffices to assume there are no elbows in $R$, which implies that the droop move $(j_1,k_1) \searrow (j_2,k_2)$ is valid. Let $B$ denote the BPD for $w$ obtained by performing this droop move on $\overline{B}_w$. Note that the tiles in $B$ and $\overline{B}_w$ only differ inside of $R$, and more specifically, we have the following reassignment of tiles on the corners of $R$:
\[
B(j_2,k_1) = \inlineare \; , \quad B(j_1,k_2) = \inlineare \; , \quad \text{and} \quad B(j_2,k_2) = \inlinejay.
\]
Consider the following two rectangles on $B$:
\[
R_1 = [j_2,j_3] \times [k_1,k_2] \subseteq B \quad \text{and} \quad R_2 = [j_1,j_2] \times [k_2,k_3] \subseteq B.
\]
Observe that $R_1 \cap R_2 = \{(j_2,k_2)\}$. Moreover, we claim that both $R_1$ and $R_2$ are active. Indeed, the only $\inlinejay$ tile in $B$ is $B(j_2,k_2)$, which is the top-right corner of $R_1$ and the bottom-left corner of $R_2$. For $R_1$, the top-left corner is $B(j_2,k_1) = \inlineare$, the bottom-left corner of $R_1$ is $B(j_3,k_1)= \inlinevwire$, and the bottom-right corner of $R_1$ is $B(j_3,k_2)=\inlinenowire$. Hence, $R_1$ is active. Likewise, for $R_2$, the top-left corner is $B(j_1,k_2) = \inlineare$, the top-right corner of $R_2$ is $B(j_1,k_3) = \inlinehwire$, and the bottom-right corner is $B(j_2,k_3) = \inlinenowire$. Hence, $R_2$ is active. Therefore, by Corollary~\ref{cor:main}, we have that $H_w$ contains a 4-cycle as advertised.
\end{proof}

\begin{proposition}
\label{p:avoid-24153}
	If $w$ contains a 24153 pattern, then $H_w$ contains a $4$-cycle.
\end{proposition}
\begin{proof}
Let $k_1 < \ldots < k_5$ be an occurrence of a 24153 pattern in $w$, so that
\[
w(k_2) < w(k_4) < w(k_1) < w(k_5) < w(k_3).
\]
Consider $\overline{B}_w$, the Rothe BPD for $w$. Use $j_1 < \ldots < j_5$ to label the five rows of $\overline{B}_w$ corresponding to the unique $\inlineare$ tiles of this pattern. Observe that $k_1 < \ldots < k_5$ label the five columns of $\overline{B}_w$ corresponding to the unique $\inlineare$ tiles of this pattern. Consider the following two rectangles on $\overline{B}_w$:
\[
R_1 = [j_1,j_2] \times [k_2,k_3] \subseteq \overline{B}_w \quad \text{and} \quad R_2 = [j_3,j_4] \times [k_1,k_3] \subseteq \overline{B}_w.
\]
Observe that $R_1$ and $R_2$ are disjoint since $j_2 < j_3$. Moreover, we claim that both $R_1$ and $R_2$ are active. Indeed, there are no $\inlinejay$ tiles in either $R_1$ or $R_2$, since they are rectangles on $\overline{B}_w$. For $R_1$, the top-left corner is $\overline{B}_w(j_1,k_2) = \inlineare$, the top-right corner of $R_1$ is $\overline{B}_w(j_1,k_3)= \inlinehwire$, the bottom-left corner of $R_1$ is $\overline{B}_w(j_2,k_2) = \inlinevwire$, and the bottom-right corner is $\overline{B}_w(j_2,k_3)=\inlinenowire$. Hence, $R_1$ is active. Likewise, for $R_2$, the top-left corner is $\overline{B}_w(j_3,k_1) = \inlineare$, the top-right corner is $\overline{B}_w(j_3,k_3) = \inlinehwire$, the bottom-left corner is $\overline{B}_w(j_4,k_1) = \inlinevwire$, and the bottom-right corner is $B(j_4,k_3) = \inlinenowire$. Hence, $R_2$ is active. Therefore, by Corollary~\ref{cor:main}, we have that $H_w$ contains a 4-cycle as advertised.
\end{proof}

\begin{proposition}
\label{p:avoid-31524}
	If $w$ contains a 31524 pattern, then $H_w$ contains a $4$-cycle.
\end{proposition}
\begin{proof}
Let $k_1 < \ldots < k_5$ be an occurrence of a 31524 pattern in $w$, so that
\[
w(k_3) < w(k_1) < w(k_5) < w(k_2) < w(k_4).
\]
Consider $\overline{B}_w$, the Rothe BPD for $w$. Use $j_1 < \ldots < j_5$ to label the five rows of $\overline{B}_w$ corresponding to the unique $\inlineare$ tiles of this pattern. Observe that $k_1 < \ldots < k_5$ label the five columns of $\overline{B}_w$ corresponding to the unique $\inlineare$ tiles of this pattern. Consider the following two rectangles on $\overline{B}_w$:
\[
R_1 = [j_1,j_3] \times [k_3,k_4] \subseteq \overline{B}_w \quad \text{and} \quad R_2 = [j_2,j_3] \times [k_1,k_2] \subseteq \overline{B}_w.
\]
Observe that $R_1$ and $R_2$ are disjoint since $k_2 < k_3$. Moreover, we claim that both $R_1$ and $R_2$ are active. Indeed, there are no $\inlinejay$ tiles in either $R_1$ or $R_2$, since they are rectangles on $\overline{B}_w$. For $R_1$, the top-left corner is $\overline{B}_w(j_1,k_3) = \inlineare$, the top-right corner of $R_1$ is $\overline{B}_w(j_1,k_4)= \inlinehwire$, the bottom-left corner of $R_1$ is $\overline{B}_w(j_3,k_3) = \inlinevwire$, and the bottom-right corner is $\overline{B}_w(j_3,k_4)=\inlinenowire$. Hence, $R_1$ is active. Likewise, for $R_2$, the top-left corner is $\overline{B}_w(j_2,k_1) = \inlineare$, the top-right corner is $\overline{B}_w(j_2,k_2) = \inlinehwire$, the bottom-left corner is $\overline{B}_w(j_3,k_1) = \inlinevwire$, and the bottom-right corner is $B(j_3,k_2) = \inlinenowire$. Hence, $R_2$ is active. Therefore, by Corollary~\ref{cor:main}, we have that $H_w$ contains a 4-cycle as advertised.\end{proof}

\begin{proposition}
\label{p:avoid-35142}
	If $w$ contains a 35142 pattern, then $H_w$ contains a $4$-cycle.
\end{proposition}
\begin{proof}
Let $k_1 < \ldots < k_5$ be an occurrence of a 35142 pattern in $w$, so that
\[
w(k_3) < w(k_5) < w(k_1) < w(k_4) < w(k_2).
\]
Consider $\overline{B}_w$, the Rothe BPD for $w$. Use $j_1 < \ldots < j_5$ to label the five rows of $\overline{B}_w$ corresponding to the unique $\inlineare$ tiles of this pattern. Observe that $k_1 < \ldots < k_5$ label the five columns of $\overline{B}_w$ corresponding to the unique $\inlineare$ tiles of this pattern. Consider the following two rectangles on $\overline{B}_w$:
\[
R_1 = [j_1,j_2] \times [k_3,k_4] \subseteq \overline{B}_w \quad \text{and} \quad R_2 = [j_3,j_4] \times [k_1,k_2] \subseteq \overline{B}_w.
\]
Observe that $R_1$ and $R_2$ are disjoint since $k_2 < k_3$. Moreover, we claim that both $R_1$ and $R_2$ are active. Indeed, there are no $\inlinejay$ tiles in either $R_1$ or $R_2$, since they are rectangles on $\overline{B}_w$. For $R_1$, the top-left corner is $\overline{B}_w(j_1,k_3) = \inlineare$, the top-right corner of $R_1$ is $\overline{B}_w(j_1,k_4)= \inlinehwire$, the bottom-left corner of $R_1$ is $\overline{B}_w(j_2,k_3) = \inlinevwire$, and the bottom-right corner is $\overline{B}_w(j_2,k_4)=\inlinenowire$. Hence, $R_1$ is active. Likewise, for $R_2$, the top-left corner is $\overline{B}_w(j_3,k_1) = \inlineare$, the top-right corner is $\overline{B}_w(j_3,k_2) = \inlinehwire$, the bottom-left corner is $\overline{B}_w(j_4,k_1) = \inlinevwire$, and the bottom-right corner is $B(j_4,k_2) = \inlinenowire$. Hence, $R_2$ is active. Therefore, by Corollary~\ref{cor:main}, we have that $H_w$ contains a 4-cycle as advertised.
\end{proof}

\begin{proposition}
\label{p:avoid-214563}
	If $w$ contains a 214563 pattern, then $H_w$ contains a $4$-cycle.
\end{proposition}
\begin{proof}
Let $k_1 < \ldots < k_6$ be an occurrence of a 214563 pattern in $w$, so that
\[
w(k_2) < w(k_1) < w(k_4) < w(k_5) < w(k_6) < w(k_3).
\]
Consider $\overline{B}_w$, the Rothe BPD for $w$. Use $j_1 < \ldots < j_6$ to label the six rows of $\overline{B}_w$ corresponding to the unique $\inlineare$ tiles of this pattern. Observe that $k_1 < \ldots < k_6$ label the six columns of $\overline{B}_w$ corresponding to the unique $\inlineare$ tiles of this pattern. Consider the rectangle
\[
R = [j_1,j_4] \times [k_2,k_3] \subseteq \overline{B}_w.
\]
In $R$, the top-left corner is $\overline{B}_w(j_1,k_2) = \inlineare$ and the bottom-right corner is $\overline{B}_w(j_4,k_3) = \inlinenowire$ from the inversion $(k_3,k_4)$. Moreover, there are no $\inlinejay$ tiles in $R$ since it is a rectangle on $\overline{B}_w$. If there is a $\inlineare$ tile in $R - \{(j_1,k_2)\}$, it would be in some row $j_1'$ with $j_1 < j_1' < j_4$ and some column $k_2'$ with $k_2 < k_2' < k_3$. If $j_1' < j_2$, we could select $k_1 \; k_2' \; k_3 \; k_4 \; k_5 \; k_6$ as our occurrence of a 214563 pattern instead. Note that $j_1'$ cannot equal $j_2$. If $j_2 < j_1' < j_4$, then $k_1 \; k_2' \; k_3 \; k_5 \; k_6$ is an occurrence of a 12453 pattern, in which case $H_w$ has a 4-cycle by Proposition~\ref{p:avoid-12453}. Hence, it suffices to assume there are no elbows in $R$, which implies that the droop move $(j_1,k_2) \searrow (j_4,k_3)$ is valid. Let $B$ denote the BPD for $w$ obtained by performing this droop move on $\overline{B}_w$. Note that the tiles in $B$ and $\overline{B}_w$ only differ inside of $R$, and more specifically, we have the following reassignment of tiles on the corners of $R$:
\[
B(j_4,k_2) = \inlineare \; , \quad B(j_1,k_3) = \inlineare \; , \quad \text{and} \quad B(j_4,k_3) = \inlinejay.
\]
Further, we have that $B(j_3,k_2) = \inlinenowire$ since $\overline{B}_w(j_3,k_2) = \inlinevwire$, and $B(j_2,k_2)=\inlinehwire$ since $\overline{B}_w(j_2,j_2) = \inlinecross$. Consider the following two rectangles on $B$:
\[
R_1 = [j_2,j_3] \times [k_1,k_2] \subseteq B \quad \text{and} \quad R_2 = [j_4,j_5] \times [k_2,k_3] \subseteq B.
\]
Observe that $R_1 \cap R_2$ are disjoint since $j_3 < j_4$. Moreover, we claim that both $R_1$ and $R_2$ are active. Indeed, the only $\inlinejay$ tile in $B$ is $B(j_5,k_3)$, which is the top-right corner of $R_2$. For $R_1$, the top-left corner is $B(j_2,k_1) = \inlineare$, the bottom-left corner of $R_1$ is $B(j_3,k_1)= \inlinevwire$, and the bottom-right corner of $R_1$ is $B(j_3,k_2)=\inlinenowire$. Hence, $R_1$ is active. Likewise, for $R_2$, the top-left corner is $B(j_4,k_2) = \inlineare$, the bottom-left corner of $R_2$ is $B(j_5,k_2) = \inlinevwire$, and the bottom-right corner is $B(j_5,k_3) = \inlinenowire$ from the inversion $(k_3,k_6)$. Hence, $R_2$ is active. Therefore, by Corollary~\ref{cor:main}, we have that $H_w$ contains a 4-cycle as advertised.
\end{proof}

\begin{proposition}
\label{p:avoid-216345}
	If $w$ contains a 216345 pattern, then $H_w$ contains a $4$-cycle.
\end{proposition}
\begin{proof}
Let $k_1 < \ldots < k_6$ be an occurrence of a 216345 pattern in $w$, so that
\[
w(k_2) < w(k_1) < w(k_6) < w(k_3) < w(k_4) < w(k_5).
\]
Consider $\overline{B}_w$, the Rothe BPD for $w$. Use $j_1 < \ldots < j_6$ to label the six rows of $\overline{B}_w$ corresponding to the unique $\inlineare$ tiles of this pattern. Observe that $k_1 < \ldots < k_6$ label the six columns of $\overline{B}_w$ corresponding to the unique $\inlineare$ tiles of this pattern. Consider the rectangle
\[
R = [j_2,j_3] \times [k_1,k_4] \subseteq \overline{B}_w.
\]
In $R$, the top-left corner is $\overline{B}_w(j_2,k_1) = \inlineare$ and the bottom-right corner is $\overline{B}_w(j_3,k_4) = \inlinenowire$ from the inversion $(k_4,k_6)$. Moreover, there are no $\inlinejay$ tiles in $R$ since it is a rectangle on $\overline{B}_w$. If there is a $\inlineare$ tile in $R - \{(j_2,k_1)\}$, it would be in some row $j_2'$ with $j_2 < j_2' < j_3$ and some column $k_1'$ with $k_1 < k_1' < k_4$. If $k_1' < k_2$, we could select $k_1' \; k_2 \; k_3 \; k_4 \; k_5 \; k_6$ as our occurrence of a 216345 pattern instead. Note that $k_1'$ cannot equal $k_2$. If $k_2 < k_1' < k_4$, then $k_2 \; k_1' \; k_4 \; k_5 \; k_6$ is an occurrence of a 12534 pattern, in which case $H_w$ has a 4-cycle by Proposition~\ref{p:avoid-12534}. Hence, it suffices to assume there are no elbows in $R$, which implies that the droop move $(j_2,k_1) \searrow (j_3,k_4)$ is valid. Let $B$ denote the BPD for $w$ obtained by performing this droop move on $\overline{B}_w$. Note that the tiles in $B$ and $\overline{B}_w$ only differ inside of $R$, and more specifically, we have the following reassignment of tiles on the corners of $R$:
\[
B(j_3,k_1) = \inlineare \; , \quad B(j_2,k_4) = \inlineare \; , \quad \text{and} \quad B(j_3,k_4) = \inlinejay.
\]
Further, we have that $B(j_2,k_3) = \inlinenowire$ since $\overline{B}_w(j_2,k_3) = \inlinehwire$, and $B(j_2,k_2)=\inlinehwire$ since $\overline{B}_w(j_2,j_2) = \inlinecross$. Consider the following two rectangles on $B$:
\[
R_1 = [j_1,j_2] \times [k_2,k_3] \subseteq B \quad \text{and} \quad R_2 = [j_2,j_3] \times [k_4,k_5] \subseteq B.
\]
Observe that $R_1 \cap R_2$ are disjoint since $k_3 < k_4$. Moreover, we claim that both $R_1$ and $R_2$ are active. Indeed, the only $\inlinejay$ tile in $B$ is $B(j_3,k_4)$, which is the bottom-left corner of $R_2$. For $R_1$, the top-left corner is $B(j_1,k_2) = \inlineare$, the top-right corner of $R_1$ is $B(j_1,k_3)= \inlinehwire$, the bottom-left corner is $B(j_2,k_2) = \inlinehwire$, and the bottom-right corner of $R_1$ is $B(j_3,k_2)=\inlinenowire$. Hence, $R_1$ is active. Likewise, for $R_2$, the top-left corner is $B(j_2,k_4) = \inlineare$, the top-right corner is $B(j_2,k_5) = \inlinehwire$, the bottom-left corner is $B(j_3,k_4) = \inlinejay$, and the bottom-right corner is $B(j_3,k_5) = \inlinenowire$ from the inversion $(k_5,k_6)$. Hence, $R_2$ is active. Therefore, by Corollary~\ref{cor:main}, we have that $H_w$ contains a 4-cycle as advertised.
\end{proof}
\begin{proposition}
\label{p:avoid-354162}
	If $w$ contains a 354162 pattern, then $H_w$ contains a $4$-cycle.
\end{proposition}
\begin{proof}
Let $k_1 < \ldots < k_6$ be an occurrence of a 354162 pattern in $w$, so that
\[
w(k_3) < w(k_5) < w(k_4) < w(k_1) < w(k_6) < w(k_2).
\]
Consider $\overline{B}_w$, the Rothe BPD for $w$. Use $j_1 < \ldots < j_6$ to label the six rows of $\overline{B}_w$ corresponding to the unique $\inlineare$ tiles of this pattern. Observe that $k_1 < \ldots < k_6$ label the six columns of $\overline{B}_w$ corresponding to the unique $\inlineare$ tiles of this pattern. Consider the following two rectangles on $\overline{B}_w$:
\[
R_1 = [j_1,j_2] \times [k_3,k_4] \subseteq \overline{B}_w \quad \text{and} \quad R_2 = [j_4,j_5] \times [k_1,k_2] \subseteq \overline{B}_w.
\]
Observe that $R_1$ and $R_2$ are disjoint since $k_2 < k_3$. Moreover, we claim that both $R_1$ and $R_2$ are active. Indeed, there are no $\inlinejay$ tiles in either $R_1$ or $R_2$, since they are rectangles on $\overline{B}_w$. For $R_1$, the top-left corner is $\overline{B}_w(j_1,k_3) = \inlineare$, the top-right corner of $R_1$ is $\overline{B}_w(j_1,k_4)= \inlinehwire$, the bottom-left corner of $R_1$ is $\overline{B}_w(j_2,k_3) = \inlinevwire$, and the bottom-right corner is $\overline{B}_w(j_2,k_4)=\inlinenowire$. Hence, $R_1$ is active. Likewise, for $R_2$, the top-left corner is $\overline{B}_w(j_4,k_1) = \inlineare$, the top-right corner is $\overline{B}_w(j_4,k_2) = \inlinehwire$, the bottom-left corner is $\overline{B}_w(j_5,k_1) = \inlinevwire$, and the bottom-right corner is $B(j_5,k_2) = \inlinenowire$. Hence, $R_2$ is active. Therefore, by Corollary~\ref{cor:main}, we have that $H_w$ contains a 4-cycle as advertised.
\end{proof}

\begin{proposition}
\label{p:avoid-461325}
	If $w$ contains a 461325 pattern, then $H_w$ contains a $4$-cycle.
\end{proposition}
\begin{proof}
Let $k_1 < \ldots < k_6$ be an occurrence of a 461325 pattern in $w$, so that
\[
w(k_4) < w(k_6) < w(k_1) < w(k_3) < w(k_2) < w(k_5).
\]
Consider $\overline{B}_w$, the Rothe BPD for $w$. Use $j_1 < \ldots < j_6$ to label the six rows of $\overline{B}_w$ corresponding to the unique $\inlineare$ tiles of this pattern. Observe that $k_1 < \ldots < k_6$ label the six columns of $\overline{B}_w$ corresponding to the unique $\inlineare$ tiles of this pattern. Consider the following two rectangles on $\overline{B}_w$:
\[
R_1 = [j_1,j_2] \times [k_4,k_5] \subseteq \overline{B}_w \quad \text{and} \quad R_2 = [j_3,j_4] \times [k_1,k_2] \subseteq \overline{B}_w.
\]
Observe that $R_1$ and $R_2$ are disjoint since $k_2 < k_4$. Moreover, we claim that both $R_1$ and $R_2$ are active. Indeed, there are no $\inlinejay$ tiles in either $R_1$ or $R_2$, since they are rectangles on $\overline{B}_w$. For $R_1$, the top-left corner is $\overline{B}_w(j_1,k_4) = \inlineare$, the top-right corner of $R_1$ is $\overline{B}_w(j_1,k_5)= \inlinehwire$, the bottom-left corner of $R_1$ is $\overline{B}_w(j_2,k_4) = \inlinevwire$, and the bottom-right corner is $\overline{B}_w(j_2,k_5)=\inlinenowire$. Hence, $R_1$ is active. Likewise, for $R_2$, the top-left corner is $\overline{B}_w(j_3,k_1) = \inlineare$, the top-right corner is $\overline{B}_w(j_3,k_2) = \inlinehwire$, the bottom-left corner is $\overline{B}_w(j_4,k_1) = \inlinevwire$, and the bottom-right corner is $B(j_4,k_2) = \inlinenowire$. Hence, $R_2$ is active. Therefore, by Corollary~\ref{cor:main}, we have that $H_w$ contains a 4-cycle as advertised.
\end{proof}
\begin{proposition}
\label{p:avoid-465132}
	If $w$ contains a 465132 pattern, then $H_w$ contains a $4$-cycle.
\end{proposition}
\begin{proof}
Let $k_1 < \ldots < k_6$ be an occurrence of a 465132 pattern in $w$, so that
\[
w(k_4) < w(k_6) < w(k_5) < w(k_1) < w(k_3) < w(k_2).
\]
Consider $\overline{B}_w$, the Rothe BPD for $w$. Use $j_1 < \ldots < j_6$ to label the six rows of $\overline{B}_w$ corresponding to the unique $\inlineare$ tiles of this pattern. Observe that $k_1 < \ldots < k_6$ label the six columns of $\overline{B}_w$ corresponding to the unique $\inlineare$ tiles of this pattern. Consider the following two rectangles on $\overline{B}_w$:
\[
R_1 = [j_1,j_2] \times [k_4,k_5] \subseteq \overline{B}_w \quad \text{and} \quad R_2 = [j_4,j_5] \times [k_1,k_2] \subseteq \overline{B}_w.
\]
Observe that $R_1$ and $R_2$ are disjoint since $k_2 < k_4$. Moreover, we claim that both $R_1$ and $R_2$ are active. Indeed, there are no $\inlinejay$ tiles in either $R_1$ or $R_2$, since they are rectangles on $\overline{B}_w$. For $R_1$, the top-left corner is $\overline{B}_w(j_1,k_4) = \inlineare$, the top-right corner of $R_1$ is $\overline{B}_w(j_1,k_5)= \inlinehwire$, the bottom-left corner of $R_1$ is $\overline{B}_w(j_2,k_4) = \inlinevwire$, and the bottom-right corner is $\overline{B}_w(j_2,k_5)=\inlinenowire$. Hence, $R_1$ is active. Likewise, for $R_2$, the top-left corner is $\overline{B}_w(j_4,k_1) = \inlineare$, the top-right corner is $\overline{B}_w(j_4,k_2) = \inlinehwire$, the bottom-left corner is $\overline{B}_w(j_5,k_1) = \inlinevwire$, and the bottom-right corner is $B(j_5,k_2) = \inlinenowire$. Hence, $R_2$ is active. Therefore, by Corollary~\ref{cor:main}, we have that $H_w$ contains a 4-cycle as advertised.
\end{proof}
\section{Future Directions}
\label{s:conc}
Here we mention some related work and briefly mention future directions for research.
\subsection{Connection to other works}
Reduced bumpless pipe dreams were introduced to study Schubert polynomials, but these were not the first such objects. One of the earliest combinatorial formulas for Schubert polynomials uses the {\em reduced pipe dreams} of~\cite{bergeron1993rc}. These also resemble tilings, but they follow different rules and use $\inlinecross$ and $\inlinebump$ tiles instead. The latter tile is called a ``bump'' tile, and that is where ``bumpless'' pipe dreams get their name. Similar to before, reduced pipe dreams of size $n$ can be partitioned into a collection indexed by permutations $w \in S_n$. Denoting the reduced pipe dreams associated to $w$ by $\text{PD}(w)$, Rubey~\cite{RUBEY} generalized the ``chute moves'' from ~\cite{bergeron1993rc} to define a partial order $\leq_{\text{chute}}$ on $\text{PD}(w)$, where $P_1 \leq_{\text{chute}} P_2$ whenever $P_2$ is obtained from $P_1$ as a chute move. \\

In~\cite{defant}, the authors prove Rubey's conjecture: the chute move poset associated to $w$ is a lattice for any permutation $w$. Since there is a known bijection $\varphi$ between $\text{PD}(w)$ and $\text{BPD}(w)$ for any $w$ (see~\cite{gao2021canonical}), one could wonder if $\varphi$ extends to a poset isomorphism. Alas, this is not the case, as can be seen even for just $w = 1432 \in S_4$. Droop moves and chute moves are inherently different. Whereas the chute move poset associated to $w$ is a lattice, the min-droop move poset $P_w$ need not have a unique minimal element. \\

Perhaps there is an analogous bumpless pipe dream version of what~\cite{defant} did for pipe dreams, so we include this as Problem~\ref{prob:poset}.
\subsection{Sufficient Conditions for Trees}
In Theorem~\ref{t:pattern-char} we show that avoiding each of the 15 patterns listed there is necessary in order for $H_w$ to be a tree. It is natural to ask whether or not avoiding these 15 patterns is sufficient. We have checked all permutations $w$ in $S_n$ for $n \leq 9$ and found that avoidance of these fifteen was sufficient for these $H_w$ graphs to be trees, but we have been unable to prove this in general so far. Hence, we list this as Problem~\ref{prob:suff}. \\

\subsection{Other Properties}
What other relationships between $w$ and $H_w$ can be proven? 
\begin{problem}
\label{prob:gao}
If $w$ contains $u$ as a pattern, does $H_w$ contains a subgraph isomorphic to $H_u$? Moreover, does having two distinct occurrences of $u$ in $w$ lead to two distinct $H_u$ subgraphs in $H_w$?
\end{problem}
Our proofs in Section~\ref{s:nec-to-avoid} above give a weaker form of this question. Progress on Problem~\ref{prob:gao} could potentially be used to attack the Gao conjecture~\cite{gao2021} or at least improve the known lower bounds on $|\text{BPD}(w)|$. 
\begin{problem}
What properties of $H_w$ (other than acyclicity) can be determined from $w$? For example, can the diameter (number of edges in the longest path) of $H_w$ be determined from $w$? Is the graph $H_w$ bipartite for any $w$?
\end{problem}

\newpage
\bibliographystyle{alpha}

\begin{thebibliography}{10}
\bibitem[AFDMT25]{defant} I. Axelrod-Freed, Colin Defant, Hanna Mularczyk, Son Nguyen, and Katerine Tung. Chute Move posets are lattices. {\em arXiv: 2507.13214}, 2025.

\bibitem[BB93]{bergeron1993rc} N. Bergeron and S. Billey. RC-graphs and Schubert polynomials. {\em Experiment. Math.}, 2(4): 257--269, 1993.

\bibitem[Gao21]{gao2021} Y. Gao. Principal specializations of Schubert polynomials and pattern containment. {\em European J. Combin.}, Volume 94, Paper No. 103291, 2021.

\bibitem[GH23]{gao2021canonical} Y. Gao and D. Huang. The canonical bijection between pipe dreams and bumpless pipe dreams. {\em International Mathematics Research Notices}, Volume 2023, Issue 21, pp. 18629--18663, 2023.

\bibitem[HP22]{huang2022bumpless} D. Huang and P. Pylyavskyy. Bumpless pipe dream RSK, growth diagrams, and Schubert structure coefficients. {\em arXiv: 2206.1435}, 2022.

\bibitem[LLS21]{lam2021back} T. Lam, S. Lee, and M. Shimozono. Back stable Schubert calculus. {\em Compos. Math.}, 157(5): 883--962, 2021.

\bibitem[Rub12]{RUBEY} M. Rubey. Maximal 0--1-fillings of moon polyominoes with restricted chain lengths and RC-graphs. {\em Advances in Applied Mathematics}, 48(2): 290--305, 2012.

\bibitem[Wei25]{weigandt2025} Anna Weigandt. Changing bases with pipe dream combinatorics. {\em arXiv: 2506.07306}, 2025.
\end{thebibliography}

\end{document}